\documentclass[submission,copyright,creativecommons]{eptcs}
\providecommand{\event}{The 13th International Conference on \\Quantum Physics and Logic (QPL) 2016}

\usepackage[cmex10]{amsmath}
\usepackage{amsthm,amsfonts,amssymb,mathtools,mathdots}

\usepackage{enumitem, framed}
\PassOptionsToPackage{svgnames}{xcolor}
\usepackage{tikz}
\usepackage{stmaryrd}
\usepackage[title,titletoc,toc]{appendix}

\newtheorem{thm}{Theorem}[section]

\theoremstyle{definition}
\newtheorem{dfn}[thm]{Definition}
\theoremstyle{remark}
\newtheorem{remark}[thm]{Remark}
\newtheorem{exm}[thm]{Example}

\newcommand\cp[1]{\,\#_{#1}\,}
\newcommand\bord[2]{\partial_{#1}^{#2}}

\usetikzlibrary{arrows, decorations.markings, shapes.geometric, decorations.pathmorphing, decorations.pathreplacing, intersections, patterns}

\pgfdeclarelayer{bg}
\pgfdeclarelayer{mid}
\pgfdeclarelayer{fg}
\pgfsetlayers{bg,mid,fg,main}

\tikzset{
	0c/.style={circle, draw, fill, inner sep=1.5pt},
	1c/.style={->, thick, shorten <=2pt, shorten >=2pt},
	1cinc/.style={right hook->, thick, shorten <=2pt, shorten >=2pt},
	2c/.style={double, thick, shorten <=10pt, shorten >=10pt, decoration={markings,mark=at position -8pt with {\arrow[scale=1.75]{>}}}, preaction={decorate}},
	3c1/.style={thick, double, double distance=3pt, shorten <=9pt, shorten >=11pt},
    	3c2/.style={thick, shorten <=9pt, shorten >=10pt},
	3c3/.style={shorten <=9pt, shorten >=10pt, decoration={markings,mark=at position -8pt with {\arrow[scale=3]{>}}},preaction={decorate}},
	4c1/.style={thick, double, double distance=4pt, shorten <=1pt, shorten >=2.75pt},
	4c2/.style={thick, double, double distance=1pt, shorten <=1pt, shorten >=1.25pt, decoration={markings,mark=at position -.05pt with {\arrow[scale=3,ultra thin]{>}}},preaction={decorate}},
	string/.style={scale=1},
	edge/.style={line width=1pt, color=black},
	edgedot/.style={dotted, line width=0.7pt, color=white},
	edgedotdark/.style={dotted, line width=0.7pt, color=black},
	dot/.style={circle, draw=black, line width=1pt, fill=white, inner sep=1.7pt},
	dotdark/.style={circle, line width=1pt, fill=black, inner sep=1.7pt},
	every node/.style={scale=.8},
}

\title{A Topological Perspective on \\Interacting Algebraic Theories}

\author{Amar Hadzihasanovic
\institute{University of Oxford}\email{amarh@cs.ox.ac.uk}}

\def\titlerunning{A Topological Perspective on Interacting Algebraic Theories}
\def\authorrunning{A. Hadzihasanovic}

\begin{document}
\maketitle

\begin{abstract}
Techniques from higher categories and higher-dimensional rewriting are becoming increasingly important for understanding the finer, computational properties of higher algebraic theories that arise, among other fields, in quantum computation. These theories have often the property of containing simpler sub-theories, whose interaction is regulated in a limited number of ways, which reveals a topological substrate when pictured by string diagrams. By exploring the double nature of computads as presentations of higher algebraic theories, and combinatorial descriptions of ``directed spaces'', we develop a basic language of directed topology for the compositional study of algebraic theories. We present constructions of computads, all with clear analogues in standard topology, that capture in great generality such notions as homomorphisms and actions, and the interactions of monoids and comonoids that lead to the theory of Frobenius algebras and of bialgebras. After a number of examples, we describe how a fragment of the ZX calculus can be reconstructed in this framework.
\end{abstract}

\section{Introduction}
A traditional presentation of an algebraic theory consists of a number of generating operations, together with a number of equations that they satisfy. If we are concerned with computational aspects of the presentation --- looking for normalisation procedures, for instance --- it is commonplace to replace equations with directed \emph{rewrite rules}. Then, in the analysis of critical pairs and confluences of a rewrite system, we are led to consider relations between different sequences of rewrites, which can in turn be relaxed to ``rewrites of rewrites'', and so on, leading into \emph{higher-dimensional rewriting theory} \cite{mimram2014towards}. From this perspective, the dichotomy between generators and relations in a presentation is resolved: they both become generators of a higher-dimensional algebraic theory, only differing in dimension.

The natural setting for higher-dimensional rewriting is \emph{higher category theory}, where, besides the objects (0-cells) and morphisms (1-cells) of basic category theory, there can also be $n$-cells between $(n-1)$-cells, for any $n > 0$. The use of terminology borrowed from topology is not coincidental: there is a sense in which the ``directed $n$-cells'' of higher categories behave like topological $n$-cells. This is exemplified by the successful application of methods from homology theory in the study of rewriting systems, based on this analogy \cite{lafont2007algebra, lafont2009polygraphic}; however, it is perhaps best pictured through the use of \emph{string diagrams} \cite{selinger2011survey, hinze2016equational} (or, more recently, surface diagrams \cite{dunn2016surface}) for reasoning about higher categories.

In one especially relevant application, string diagrams have emerged as a strong contender for a high-level, native syntax for quantum programming \cite{zeng2015abstract, coecke2015picturing}, whose highly symmetrical semantics --- in pure, finite-dimensional quantum theory, all processes are reversible, and inputs can be turned into outputs and vice versa \cite{selinger2011finite} --- require a quite unusual amount of interplay between algebraic and coalgebraic structures.

Better understanding the computational properties of theories such as the ZX calculus \cite{coecke2008interacting, backens2014zx} and its refinement, the ZW calculus \cite{hadzihasanovic2015diagrammatic}, is pivotal in making them viable for the efficient design of quantum algorithms and protocols. Although these theories include, as a whole, a relatively large number of axioms, they contain a number of simpler sub-theories, whose interactions are regulated by the axioms in fairly predictable ways: something is a homomorphism of something else, something is an action on something else...

Beyond these motivating examples, such a factorisation seems to be a property of many theories: as a simple case, think of $*$-monoids, which can be seen as an interaction of the theory of monoids and the theory of involutions. So we asked ourselves the question: 
\begin{itemize}
	\item Is there a way to study algebraic theories \emph{compositionally}, so that one can derive properties of the larger theory from its components, and the few ways in which they are allowed to interact?
\end{itemize}
There has already been, in fact, an attempt to develop a compositional algebra, through Lack's ``composing PROPs'' framework \cite{lack2004composing}. In this setting, a presentation of a fragment of the ZX calculus --- the theory of \emph{interacting bialgebras} --- was successfully constructed from the theories of monoids and of comonoids \cite{bonchi2014interacting}. There are, nevertheless, two downsides to this approach, relative to our objectives.

Firstly, composition relies on the choice of a ``distributive law'', which conceptually amounts to stating what the normal form for operations of the resulting theory should be. Thus, we can already tell what the resulting theory will \emph{globally} look like, which subtracts something from its heuristic value, especially when we only have an algebraic presentation at hand. In fact, ``composing PROPs'' is mostly useful to derive axioms when concrete models of the component theories are available, also suggesting a ``concrete'' way of composing them.

Secondly, it is a \emph{flat} composition, in that it works in a strictly 2-categorical framework, and fails to account for any of the topological properties of the interactions. For instance, two specular distributive laws for monoids and comonoids lead, respectively, to the theory of special Frobenius algebras, and to the theory of bialgebras. Both the theory of monoids and the theory of comonoids are \emph{planar} --- none of the axioms require the swapping of inputs or outputs of operations --- and so is the theory of Frobenius algebras; hence, the interaction leading to Frobenius algebras is not supposed to change the dimension of generators.

On the other hand, the bialgebra law --- a part of the theory of bialgebras --- is \emph{not} planar, and is in fact best represented by a string diagram in 3 dimensions, where the monoid part and the comonoid part belong to different, orthogonal planes, as in the following picture.
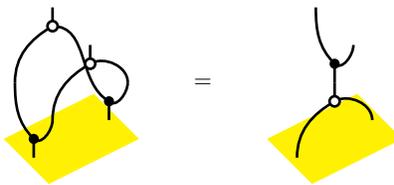
\begin{equation*}
\begin{tikzpicture}[scale=0.5] 
\begin{scope}[shift={(-3.5,0)}] 
\begin{pgfonlayer}{bg}
	\path[fill, color=yellow!100] (-1.8,0) -- (0.6,1.2) -- (1.8,0) -- (-0.6,-1.2) -- cycle;
\end{pgfonlayer}
\begin{scope}[shift={(-1,-0.5)}] 
\begin{pgfonlayer}{mid}
	\draw[edge] (0,0) to (0,0.5);
	\draw[edge, out=-45, in=-90] (0,0.5) to (0.5,1);
	\draw[edge, out=135, in=-90] (0,0.5) to (-0.5, 2);
\end{pgfonlayer}
\begin{pgfonlayer}{fg}
	\node[dotdark] at (0,0.5) {};
\end{pgfonlayer}
\end{scope}
\begin{scope}[shift={(1,0.5)}] 
\begin{pgfonlayer}{mid}
	\draw[edge] (0,0) to (0,0.5);
	\draw[edge, out=-45, in=-90] (0,0.5) to (0.5,1);
	\draw[edge, out=150, in=0] (0,0.5) to (-1.5, 2.5);
\end{pgfonlayer}
\begin{pgfonlayer}{fg}
	\node[dotdark] at (0,0.5) {};
\end{pgfonlayer}
\end{scope}
\begin{scope}[shift={(0.5,1)}] 
\begin{pgfonlayer}{mid}
	\draw[edge] (0,1) to (0,1.5);
	\draw[edge, out=90, in=-150] (-1,-0.5) to (0,1);
	\draw[edge, out=90, in=30] (1,0.5) to (0,1);
\end{pgfonlayer}
\begin{pgfonlayer}{fg}
	\node[dot] at (0,1) {};
\end{pgfonlayer}
\end{scope}
\begin{scope}[shift={(-0.5,2)}] 
\begin{pgfonlayer}{mid}
	\draw[edge] (0,1) to (0,1.5);
	\draw[edge, out=90, in=-150] (-1,-0.5) to (0,1);
\end{pgfonlayer}
\begin{pgfonlayer}{fg}
	\node[dot] at (0,1) {};
\end{pgfonlayer}
\end{scope}
\end{scope}
\begin{scope}[shift={(3.5,0)}] 
\begin{pgfonlayer}{bg}
	\path[fill, color=yellow!100] (-1.8,0) -- (0.6,1.2) -- (1.8,0) -- (-0.6,-1.2) -- cycle;
\end{pgfonlayer}
\begin{scope} 
\begin{pgfonlayer}{mid}
	\draw[edge] (0,1) to (0,1.5);
	\draw[edge, out=90, in=-150] (-1,-0.5) to (0,1);
	\draw[edge, out=90, in=30] (1,0.5) to (0,1);
\end{pgfonlayer}
\begin{pgfonlayer}{fg}
	\node[dot] at (0,1) {};
\end{pgfonlayer}
\end{scope}
\begin{scope}[shift={(0,1.5)}] 
\begin{pgfonlayer}{mid}
	\draw[edge] (0,0) to (0,0.5);
	\draw[edge, out=-45, in=-90] (0,0.5) to (0.5,1);
	\draw[edge, out=135, in=-90] (0,0.5) to (-0.5, 2);
\end{pgfonlayer}
\begin{pgfonlayer}{fg}
	\node[dotdark] at (0,0.5) {};
\end{pgfonlayer}
\end{scope}
\end{scope}
\begin{scope}[shift={(0,1.5)}] 
	\node at (0,0) {$=$};
\end{scope}
\end{tikzpicture}
\end{equation*}
So the interaction leading to bialgebras should be of a different, dimension-increasing sort. 

In this paper, we try to lay the groundwork for an alternative approach, and make a case for the following assertions:
\begin{enumerate}
	\item that there exists a way of studying algebraic theories compositionally, with a small number of basic constructions corresponding to the most frequent interactions;
	\item that the language for compositional algebra is a kind of combinatorial \emph{directed topology}, all interactions having clear analogues in standard topology.
\end{enumerate}

In Section \ref{sec:computads}, we briefly present our technical framework of choice, the theory of \emph{computads} or polygraphs, and build a basic vocabulary of directed topology in this context. In Section \ref{sec:basic}, we use these tools to construct presentations of some basic theories, such as the theory of monoids, from even simpler ones. In Section \ref{sec:homo}, we show how two kinds of interaction capture the notion of homomorphism and of (co)action. Finally, in Section \ref{sec:smash}, we introduce the dimension-changing operation which enables us to obtain the theory of bialgebras (and of commutative monoids) from the theory of monoids. We conclude by describing a partial reconstruction ``from scratch'' of the theory of interacting bialgebras, and discuss some of the many possible further directions of this project.

\section{Computads and directed spaces} \label{sec:computads}
A \emph{computad} is, informally, a presentation of a higher category ``\emph{by generators and no relations}''. This notion, introduced by Street \cite{street1976limits} in the 2-dimensional case, was rediscovered and expanded by Burroni \cite{burroni1993higher} in the context of higher-dimensional rewriting theory, where it is known under the name of \emph{polygraph}. We will only use the original setting of computads for strict $\omega$-categories, although weak variants have also been considered \cite{batanin2002computads}.

A computad describes how to build a higher category by progressively adjoining generating cells of higher dimension, whose border is a pasting of lower-dimensional cells, and letting composition be, at each stage, the free pasting of old and new cells. 

Like Lawvere theories, PROs, and operads \cite{markl2008operads}, computads can be used to describe algebraic theories, to be internalised in arbitrary higher categories through appropriate ``semantical'' functors. The equational laws by which such theories are normally presented need to have corresponding higher cells, and are, therefore, directed by default: this makes computads a natural choice for studying computational aspects of presentations (normalisation, confluence...), and has led to a particular interest among rewriting theorists \cite{guiraud2006three, lafont2007algebra}.

The information contained in a computad is conceptually analogous to the description of a topological space as a CW complex, with \emph{directed cells} replacing the undirected, topological cells. The formalism of higher categories can be seen as an auxiliary, combinatorial tool for specifying how cells are glued together, corresponding to the gluing maps of point-set topology; in other words, a computad is a combinatorial description of a \emph{directed space}, in the spirit of Grandis \cite{grandis2009directed}. There may be ways to make an explicit connection; for now, we take this as no more than a guiding heuristics.

The simplest and more direct combinatorics of pasting are provided by globular $\omega$-categories. Since technical details are not particularly important in the remainder, we leave the presentation quite informal, and refer, for instance, to \cite{metayer2003resolutions} for more details. Preliminary definitions are in Appendix \ref{sec:globular}. 

\begin{dfn}
A \emph{computad} $X$ is a pair of an $\omega$-category $\mathcal{X}$ and a subset $|X| \subseteq \mathcal{X}$, such that, if $\mathcal{X}_n$ is the $n$-skeleton of $\mathcal{X}$, and $|X|_n := \{ \sigma \in |X| \;|\; d(\sigma) = n \}$, for $n \geq 0$:
\begin{itemize}
	\item $\mathcal{X}_0 = |X|_0$;
	\item $\mathcal{X}_n$ is obtained from $\mathcal{X}_{n-1}$ by freely adjoining the cells of $|X|_n$.
\end{itemize}

A map $f: X \to Y$ of computads is a function of sets $|f|: |X| \to |Y|$ that induces a functor of free $\omega$-categories $f: \mathcal{X} \to \mathcal{Y}$. Computads and their maps form a category $\mathbf{Cpt}$.
\end{dfn}
Thus, a computad is described by giving, for each $n$, a set $X_n$ of generating $n$-cells, and specifying their border as a formal composition of the lower-dimensional generating cells. A map of computads is a mapping of the generating cells of $X$ onto generating cells of $Y$ preserving the composition of borders.

\begin{dfn}
The \emph{disjoint union} $X \oplus Y$ of two computads $X$ and $Y$ is the computad with generating set $|X \oplus Y| := |X| + |Y|$, and 
\begin{equation*}
	\partial^\alpha_n\sigma = \begin{cases}
		\partial^\alpha_{X,n} \sigma \;, & \sigma \in |X|\;, \\
		\partial^\alpha_{Y,n} \sigma \;, & \sigma \in |Y|\;.
	\end{cases}
\end{equation*}

A \emph{subcomputad} $A \subseteq X$ is a map $A \to X$ of computads whose underlying set-function is an inclusion of sets $|A| \subseteq |X|$. Given a set $B$ of cells of $X$, we also denote by $B$ the smallest subcomputad of $X$ that contains them.

An \emph{equivalence relation} $E$ on a computad $X$ is an equivalence relation on the set $|X|$ of generating cells, whose extension to composite cells commutes with all border operators. We denote by $X/E$ the \emph{quotient} of $X$ by the equivalence relation $E$, with generating set $|X|/E$. Given a subcomputad $A \subseteq X$, we also write $X/A$ for the quotient of $X$ by the equivalence relation $|A| \times |A|$.
\end{dfn}

These operations are analogous to the corresponding operations on topological spaces, something we try to highlight through notation and terminology: we build the disjoint union of $X$ and $Y$ by separately attaching the generating cells of $X$ and of $Y$, and obtain a quotient space by identifying cells compatibly with their borders. The disjoint union is a categorical coproduct in $\mathbf{Cpt}$ \cite{makkai2001comparing}; however, the \emph{product} of topological spaces does not correspond to the categorical product of computads, but to a different monoidal structure --- the computadic version of the \emph{Crans-Gray} tensor product of $\omega$-categories \cite{crans1995pasting}.

\begin{dfn}
The \emph{tensor product} $X \otimes Y$ of two computads $X$ and $Y$ is the computad with generating set $|X \otimes Y| := |X| \times |Y|$, where 
\begin{equation*}
	|X \otimes Y|_n = \sum_{k=0}^n |X|_k \times |Y|_{n-k}\;,
\end{equation*}
and, for all generating $n$-cells $\sigma \otimes \tau$, $\sigma \in |X|_k$, $\tau \in |Y|_{n-k}$, attaching is characterised by the following condition. Let 
\begin{equation*}
	\varepsilon(n) = \begin{cases} + & \text{if $n$ is even,} \\
		- & \text{if $n$ is odd.}
	\end{cases}
\end{equation*}

By induction: $(\mathcal{X}\otimes \mathcal{Y})_0$ is just $X_0 \times Y_0$. For all $n >0$, in $(\mathcal{X} \otimes \mathcal{Y})_{n-1}$, extend the $- \otimes -$ operation to composite cells by writing $(\rho \cp{m} \rho') \otimes \pi$ for the pasting of $\rho \otimes \pi$ and $\rho' \otimes \pi$ along their shared border $\partial^+_{m-1} \rho \otimes \pi = \partial^-_{m-1}\rho' \otimes \pi$, and similarly for $\rho \otimes (\pi \cp{m} \pi')$. 

Then, $\partial^\alpha_{n-1}(\sigma \otimes \tau)$ is obtained by pasting the $(n-1)$-cells 
\begin{equation*}
	\partial^\alpha_{k-1} \sigma \otimes \tau \quad \text{and} \quad \sigma \otimes \partial^{\varepsilon(k)\alpha}_{n-k-1} \tau
\end{equation*}
along their shared border $\partial^\alpha_{k-1} \sigma \otimes \partial^{\varepsilon(k)\alpha}_{n-k-1} \tau$.
\end{dfn}

The tensor product determines a non-symmetric monoidal structure on $\mathbf{Cpt}$, with $1$, the terminal computad of a single 0-cell, as unit.

The explicit combinatorics for the border operators in $X \otimes Y$ are quite complicated; we will, however, only consider low-dimensional cases, of which some expressions are given in Appendix \ref{sec:borders}. The approach to $\omega$-categories by cubical sets with connections \cite{brown1981algebra, brown1987tensor} leads to a much sleeker definition of the Crans-Gray tensor product \cite{al2002multiple}, but it is unclear whether a cubical description of computads would lead to an overall simplification, due to the additional complications related to the handling of thin cells.

\begin{remark} That this is a valid definition can be seen as following from the results of \cite{steiner2004omega}. Our conditions determine the tensor product for the category of loop-free augmented directed complexes with unital bases, which, \emph{modulo} an adjustment of terminology, is equivalent to a subcategory of $\omega$-categories presented by computads that are loop-free in a suitable sense. These include, in particular, the $n$-globes $G^n$, generated by one $n$-cell $\top$ and one cell $k^\alpha$ for all $k < n$, $\alpha \in \{+,-\}$, such that $\bord{k}{\alpha}(\top) = k^\alpha$, and the ``walking $m$-compositions'' $G^n \cp{m} G^n$ (called $G[n]$ and $G[n;m]$ in the referenced paper). 

This tensor product determines, in turn, the Crans-Gray tensor product on the whole of $\omega\mathbf{Cat}$. Let $X$, $Y$ be computads, and $\mathcal{X} \otimes \mathcal{Y}$ be the tensor product of the $\omega$-categories they generate. For all $\sigma \in |X|_k$, $\tau \in |Y|_{n-k}$, $\sigma \otimes \tau$ corresponds to a functor $(\sigma \otimes \tau)^*: G^k \otimes G^{n-k} \to \mathcal{X} \otimes \mathcal{Y}$ sending $\top \otimes \top$ to $\sigma \otimes \tau$. 

In particular, the border of $\sigma \otimes \tau$ is the image of a cell of $G^k \otimes G^{n-k}$. But since $G^k \otimes G^{n-k}$ is freely generated by the $i^\alpha \otimes j^\beta$, and $(\sigma \otimes \tau)^*(i^\alpha \otimes j^\beta) = \partial^\alpha_i\sigma \otimes \partial^\beta_j\tau$, the border of $\sigma \otimes \tau$ is generated by the $\partial^\alpha_i\sigma \otimes \partial^\beta_j\tau$. If $\partial^\alpha_i\sigma$, $\partial^\beta_j\tau$ are generators, there is nothing to do.

Otherwise, suppose the first one is a composite cell; by the interchange law, it can be written in the form $\rho \cp{m} \rho'$ for some $m > 0$, so that no composition of dimension $d > m$ appears in $\rho, \rho'$. By a similar argument as before, using the fact that $(\rho \cp{m} \rho') \otimes \pi$ corresponds to a functor $(G^i \cp{m} G^i) \otimes G^j \to \mathcal{X} \otimes \mathcal{Y}$, we can rewrite the tensor as a composition of $\rho \otimes \pi$, $\rho' \otimes \pi$, and lower-dimensional cells. Continuing like this until all the highest-dimensional compositions are eliminated, since generators can appear only in finite number in the cells of $\mathcal{X}$ and $\mathcal{Y}$, we finally obtain an expression of the border of $\sigma \otimes \tau$ as a composition of tensors of lower-dimensional generators; so $X \otimes Y$, as defined, is a presentation of $\mathcal{X} \otimes \mathcal{Y}$.
\end{remark}

All these constructions have analogues for topological spaces; the following, however, is purely directed.
\begin{dfn}
Let $X$ be a computad, $S \subseteq \mathbb{N}$. Then, $X^{\mathrm{op}(S)}$ is the computad with the same generating set of $X$, but with the direction of $n$-cells and $n$-composition reversed for all $n \in S$; that is, letting $\sigma^{\mathrm{op}(S)} := \sigma$ for all $0$-cells of $X$, define inductively
\begin{equation*}
	\partial_{n-1}^\alpha(\sigma^{\mathrm{op}(S)}) := \begin{cases}
		(\partial_{n-1}^{-\alpha} \sigma)^{\mathrm{op}(S)}\;, & n \in S\;, \\
		(\partial_{n-1}^\alpha \sigma)^{\mathrm{op}(S)}\;, & n \not\in S\;,
	\end{cases}
\end{equation*}	
\begin{equation*}
	(\sigma \cp{n} \tau)^{\mathrm{op}(S)} := \begin{cases}
		\tau^{\mathrm{op}(S)} \cp{n} \sigma^{\mathrm{op}(S)} \;, & n \in S\;, \\
		\sigma^{\mathrm{op}(S)} \cp{n} \tau^{\mathrm{op}(S)} \;, & n \not\in S\;,
	\end{cases}
\end{equation*}	
for all numbers $n$, cells $\sigma$, $\tau$ of $\mathcal{X}_n$, and $\alpha \in \{+,-\}$. Then, $X^{\mathrm{op}(S)}$ is generated by the $\sigma^{\mathrm{op}(S)}$, for all $\sigma \in |X|$. We write $X^\mathrm{op} := X^{\mathrm{op}(\mathbb{N})}$.
\end{dfn}
In the definition of the tensor product of computads, the border of the tensor of two cells is reversed when the border of the two cells, separately, is: it follows that, for all computads $X$ and $Y$, $(X \otimes Y)^\mathrm{op} \simeq X^\mathrm{op} \otimes Y^\mathrm{op}$.

In what follows, our constructions will often result in lax versions of algebraic theories, where equalities are replaced by cells pointing in a specific direction. These may not always be the best choices for computational purposes; we will leave it implicit that one can always reverse cells of a given dimension, if needed. We will also tend to not distinguish between strict and lax versions of a theory: the theories that we consider are usually interpreted in low-dimensional categories, so we can leave it to semantical functors to strictify as necessary, turning directed cells into equalities or isomorphisms.

\section{Basic examples: Yang-Baxter, associativity, units} \label{sec:basic}

In this section, we look at some fundamental examples of interacting computads, with a focus on their algebraic interpretation. The latter is better understood when cells are visualised as \emph{string diagrams} in the style of \cite{hinze2016equational}; for the first examples, however, we will also provide the more traditional, dual presentation by pasting diagrams.

The simplest 1-dimensional computad, and the basic ingredient of many later constructions, is the directed interval $I$, which is the same as the $1$-globe $G^1$.
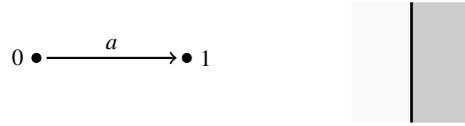
\begin{equation*}
\begin{tikzpicture}[baseline=-25pt]
\begin{pgfonlayer}{mid}
	\node[0c] (0) at (0,0) [label=left:$0$] {}; 
	\node[0c] (1) at (2,0) [label=right:$1$] {};
	\draw[1c] (0) to node[above] {$a$} (1);
\end{pgfonlayer} 
\end{tikzpicture}
\hspace{50pt}
\begin{tikzpicture}[scale=0.8]
\begin{pgfonlayer}{bg}
	\fill[color=gray!5] (-1,-1) rectangle (0,1);
	\fill[color=gray!40] (0,-1) rectangle (1,1);
\end{pgfonlayer}
\begin{pgfonlayer}{mid}
	\draw[edge] (0,-1) to (0,1);
\end{pgfonlayer}
\end{tikzpicture}
\end{equation*}

\begin{dfn}
The \emph{directed $n$-cube} is the computad $I^n$ obtained by tensoring $n$ copies of $I$.
\end{dfn}
The directed 2-cube models an ``interaction'' operation with no additional relations. We write $\sigma \tau$ for $\sigma \otimes \tau$ in the pasting diagrams.

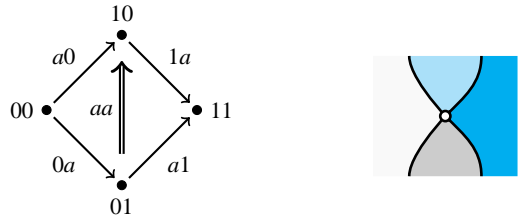
\begin{equation*}
\begin{tikzpicture}
\begin{pgfonlayer}{mid}
	\node[0c] (00) at (-1,0) [label=left:$00$] {}; 
	\node[0c] (01) at (0,-1) [label=below:$01$] {};
	\node[0c] (10) at (0,1) [label=above:$10$] {};
	\node[0c] (11) at (1,0) [label=right:$11$] {};
	\draw[1c] (00) to node[auto,swap] {$0a$} (01);
	\draw[1c] (01) to node[auto,swap] {$a1$} (11);
	\draw[1c] (00) to node[auto] {$a0$} (10);
	\draw[1c] (10) to node[auto] {$1a$} (11);
	\draw[2c] (01) to node[auto] {$aa$} (10);
\end{pgfonlayer} 
\end{tikzpicture}
\hspace{50pt}
\begin{tikzpicture}[baseline=-40pt, scale=0.8]
\begin{pgfonlayer}{bg}
	\path[fill, color=gray!5] (-1.2,-1) -- (-0.6,-1) to [out=90, in=45] (0,0) to [out=-45, in= -90] (-0.6,1) -- (-1.2,1) -- cycle;
	\path[fill, color=gray!40] (-0.6,-1) to [out=90, in=45] (0,0) to [out=-45, in= 90] (0.6,-1) -- cycle;
	\path[fill, color=cyan!30] (-0.6,1) to [out=-90, in=-45] (0,0) to [out=45, in= -90] (0.6,1) -- cycle;
	\path[fill, color=cyan!100] (1.2,-1) -- (0.6,-1) to [out=90, in=-45] (0,0) to [out=45, in= -90] (0.6,1) -- (1.2,1) -- cycle;
\end{pgfonlayer}
\begin{pgfonlayer}{mid}
	\draw[edge, out=90, in=45] (-0.6, -1) to (0,0);
	\draw[edge, out=45, in=-90] (0, 0) to (0.6,1);
	\draw[edge, out=90, in=-45] (0.6, -1) to (0,0);
	\draw[edge, out=-45, in=-90] (0, 0) to (-0.6,1);
\end{pgfonlayer}
\begin{pgfonlayer}{fg}
	\node[dot] at (0,0) {};
\end{pgfonlayer}
\end{tikzpicture}
\end{equation*}
Further on, the directed 3-cube can be seen as the presentation of a coloured, directed version of the \emph{Yang-Baxter equation}.
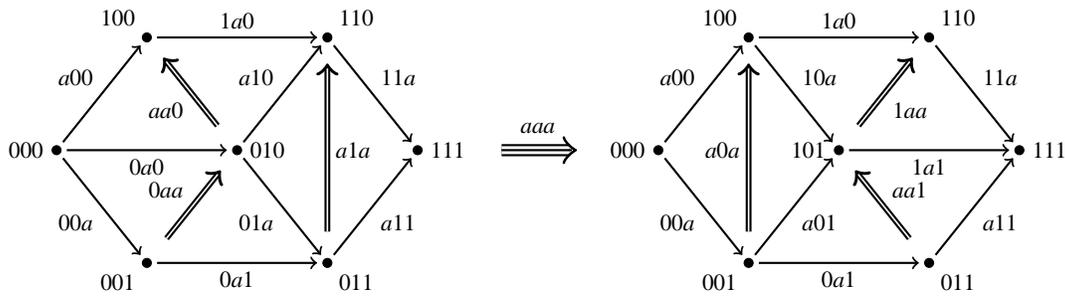
\begin{equation*}
\begin{tikzpicture}
\begin{pgfonlayer}{mid}
\begin{scope}[shift={(-4,0)}]
	\node[0c] (000) at (-2.4,0) [label=left:$000$] {}; 
	\node[0c] (001) at (-1.2,-1.5) [label=below left:$001$] {};
	\node[0c] (011) at (1.2,-1.5) [label=below right:$011$] {};
	\node[0c] (111) at (2.4,0) [label=right:$111$] {};
	\node[0c] (100) at (-1.2,1.5) [label=above left:$100$] {};
	\node[0c] (110) at (1.2,1.5) [label=above right:$110$] {};
	\node[0c] (010) at (0,0) [label=right:$010$] {};
	\draw[1c] (000) to node[auto,swap] {$00a$} (001); 
	\draw[1c] (001) to node[auto,swap] {$0a1$} (011);
	\draw[1c] (011) to node[auto,swap] {$a11$} (111);
	\draw[1c] (000) to node[auto] {$a00$} (100);
	\draw[1c] (100) to node[auto] {$1a0$} (110);
	\draw[1c] (110) to node[auto] {$11a$} (111);
	\draw[1c] (000) to node[auto,swap] {$0a0$} (010); 
	\draw[1c] (010) to node[auto,swap] {$01a$} (011);
	\draw[1c] (010) to node[auto] {$a10$} (110);
	\draw[2c] (001) to node[auto] {$0aa$} (010); 
	\draw[2c] (011) to node[auto,swap] {$a1a$} (110);
	\draw[2c] (010) to node[auto] {$aa0$} (100);
\end{scope}
\begin{scope}
	\node at (0,0.3) {$aaa$};
	\draw[3c1] (-0.8,0) to (0.8,0);
	\draw[3c2] (-0.8,0) to (0.8,0);
	\draw[3c3] (-0.8,0) to (0.8,0);
\end{scope}
\begin{scope}[shift={(4,0)}]
	\node[0c] (000) at (-2.4,0) [label=left:$000$] {}; 
	\node[0c] (001) at (-1.2,-1.5) [label=below left:$001$] {};
	\node[0c] (011) at (1.2,-1.5) [label=below right:$011$] {};
	\node[0c] (111) at (2.4,0) [label=right:$111$] {};
	\node[0c] (100) at (-1.2,1.5) [label=above left:$100$] {};
	\node[0c] (110) at (1.2,1.5) [label=above right:$110$] {};
	\node[0c] (101) at (0,0) [label=left:$101$] {};
	\draw[1c] (000) to node[auto,swap] {$00a$} (001); 
	\draw[1c] (001) to node[auto,swap] {$0a1$} (011);
	\draw[1c] (011) to node[auto,swap] {$a11$} (111);
	\draw[1c] (000) to node[auto] {$a00$} (100);
	\draw[1c] (100) to node[auto] {$1a0$} (110);
	\draw[1c] (110) to node[auto] {$11a$} (111);
	\draw[1c] (001) to node[auto,swap] {$a01$} (101); 
	\draw[1c] (101) to node[auto,swap] {$1a1$} (111);
	\draw[1c] (100) to node[auto] {$10a$} (101);
	\draw[2c] (011) to node[auto,swap] {$aa1$} (101); 
	\draw[2c] (001) to node[auto] {$a0a$} (100);
	\draw[2c] (101) to node[auto,swap] {$1aa$} (110);
\end{scope}
\end{pgfonlayer} 
\end{tikzpicture}
\end{equation*}
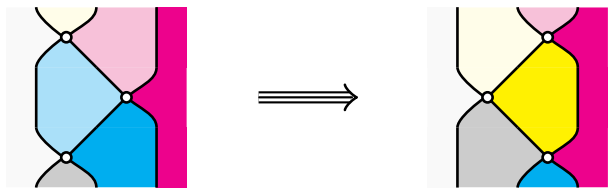
\begin{equation} \label{eq:yangbaxter}
\begin{tikzpicture}[string, yscale=0.4, xscale=0.8]
\begin{scope}[shift={(-3.5,0)}]
\begin{pgfonlayer}{bg}
\begin{scope}[shift={(-0.5,1.99)}]
	\path[fill, color=gray!5] (-1,-1) -- (-0.5,-1) to [out=90, in=60] (0,0) to [out=-60, in= -90] (-0.5,1) -- (-1,1) -- cycle;
	\path[fill, color=cyan!30] (-0.5,-1) to [out=90, in=60] (0,0) to (0.5,-1) -- cycle;
	\path[fill, color=yellow!10] (-0.5,1) to [out=-90, in=-60] (0,0) to [out=60, in= -90] (0.5,1) -- cycle;
	\path[fill, color=magenta!30] (1.5,-1) -- (0.5,-1) to (0,0) to [out=60, in= -90] (0.5,1) -- (1.5,1) -- cycle;
	\path[fill, color=magenta!100] (1.5,-1) rectangle (2,1);
\end{scope}
\begin{scope}[shift={(0.5,0)}]
	\path[fill, color=gray!5] (-2,-1) rectangle (-1.5,1);
	\path[fill, color=cyan!30] (-1.5,-1) -- (-0.5,-1) to (0,0) to (-0.5,1) -- (-1.5,1) -- cycle;
	\path[fill, color=cyan!100] (-0.5,-1) to (0,0) to [out=-60, in= 90] (0.5,-1) -- cycle;
	\path[fill, color=magenta!30] (-0.5,1) to (0,0) to [out=60, in= -90] (0.5,1) -- cycle;
	\path[fill, color=magenta!100] (1,-1) -- (0.5,-1) to [out=90, in=-60] (0,0) to [out=60, in= -90] (0.5,1) -- (1,1) -- cycle;
\end{scope}
\begin{scope}[shift={(-0.5,-1.99)}]
	\path[fill, color=gray!5] (-1,-1) -- (-0.5,-1) to [out=90, in=60] (0,0) to [out=-60, in= -90] (-0.5,1) -- (-1,1) -- cycle;
	\path[fill, color=gray!40] (-0.5,-1) to [out=90, in=60] (0,0) to [out=-60, in= 90] (0.5,-1) -- cycle;
	\path[fill, color=cyan!30] (-0.5,1) to [out=-90, in=-60] (0,0) to (0.5,1) -- cycle;
	\path[fill, color=cyan!100] (1.5,-1) -- (0.5,-1) to [out=90, in=-60] (0,0) to (0.5,1) -- (1.5,1) -- cycle;
	\path[fill, color=magenta!100] (1.5,-1) rectangle (2,1);
\end{scope}
\end{pgfonlayer}
\begin{pgfonlayer}{mid}
\begin{scope}[shift={(-0.5,2)}]
	\draw[edge, out=90, in=60] (-0.5, -1) to (0,0);
	\draw[edge, out=60, in=-90] (0, 0) to (0.5,1);
	\draw[edge] (0.5, -1) to (0,0);
	\draw[edge, out=-60, in=-90] (0, 0) to (-0.5,1);
	\draw[edge] (1.5,-1) to (1.5,1);
	\node[dot] at (0,0) {};
\end{scope}
\begin{scope}[shift={(0.5,0)}]
	\draw[edge] (-0.5, -1) to (0,0);
	\draw[edge, out=60, in=-90] (0, 0) to (0.5,1);
	\draw[edge, out=90, in=-60] (0.5, -1) to (0,0);
	\draw[edge] (0, 0) to (-0.5,1);
	\draw[edge] (-1.5,-1) to (-1.5,1);
	\node[dot] at (0,0) {};
\end{scope}
\begin{scope}[shift={(-0.5,-2)}]
	\draw[edge, out=90, in=60] (-0.5, -1) to (0,0);
	\draw[edge] (0, 0) to (0.5,1);
	\draw[edge, out=90, in=-60] (0.5, -1) to (0,0);
	\draw[edge, out=-60, in=-90] (0, 0) to (-0.5,1);
	\draw[edge] (1.5,-1) to (1.5,1);
	\node[dot] at (0,0) {};
\end{scope}
\end{pgfonlayer}
\end{scope}
\begin{scope}
	\draw[3c1] (-1.2,0) to (1.2,0);
	\draw[3c2] (-1.2,0) to (1.2,0);
	\draw[3c3] (-1.2,0) to (1.2,0);
\end{scope}
\begin{scope}[shift={(3.5,0)}]
\begin{pgfonlayer}{bg}
\begin{scope}[shift={(0.5,1.99)}]
	\path[fill, color=gray!5] (-2,-1) rectangle (-1.5,1);
	\path[fill, color=yellow!10] (-1.5,-1) -- (-0.5,-1) to (0,0) to [out=-60, in= -90] (-0.5,1) -- (-1.5,1) -- cycle;
	\path[fill, color=yellow!100] (-0.5,-1) to (0,0) to [out=-60, in= 90] (0.5,-1) -- cycle;
	\path[fill, color=magenta!30] (-0.5,1) to [out=-90, in=-60] (0,0) to [out=60, in= -90] (0.5,1) -- cycle;
	\path[fill, color=magenta!100] (1,-1) -- (0.5,-1) to [out=90, in=-60] (0,0) to [out=60, in= -90] (0.5,1) -- (1,1) -- cycle;
\end{scope}
\begin{scope}[shift={(-0.5,0)}]
	\path[fill, color=gray!5] (-1,-1) -- (-0.5,-1) to [out=90, in=60] (0,0) to [out=-60, in= -90] (-0.5,1) -- (-1,1) -- cycle;
	\path[fill, color=gray!40] (-0.5,-1) to [out=90, in=60] (0,0) to (0.5,-1) -- cycle;
	\path[fill, color=yellow!10] (-0.5,1) to [out=-90, in=-60] (0,0) to (0.5,1) -- cycle;
	\path[fill, color=yellow!100] (1.5,-1) -- (0.5,-1) to (0,0) to (0.5,1) -- (1.5,1) -- cycle;
	\path[fill, color=magenta!100] (1.5,-1) rectangle (2,1);
\end{scope}
\begin{scope}[shift={(0.5,-1.99)}]
	\path[fill, color=gray!5] (-2,-1) rectangle (-1.5,1);
	\path[fill, color=gray!40] (-1.5,-1) -- (-0.5,-1) to [out=90, in=60] (0,0) to (-0.5,1) -- (-1.5,1) -- cycle;
	\path[fill, color=cyan!100] (-0.5,-1) to [out=90, in=60] (0,0) to [out=-60, in= 90] (0.5,-1) -- cycle;
	\path[fill, color=yellow!100] (-0.5,1) to (0,0) to [out=60, in= -90] (0.5,1) -- cycle;
	\path[fill, color=magenta!100] (1,-1) -- (0.5,-1) to [out=90, in=-60] (0,0) to [out=60, in= -90] (0.5,1) -- (1,1) -- cycle;
\end{scope}
\end{pgfonlayer}
\begin{pgfonlayer}{mid}
\begin{scope}[shift={(0.5,2)}]
	\draw[edge] (-0.5, -1) to (0,0);
	\draw[edge, out=60, in=-90] (0, 0) to (0.5,1);
	\draw[edge, out=90, in=-60] (0.5, -1) to (0,0);
	\draw[edge, out=-60, in=-90] (0, 0) to (-0.5,1);
	\draw[edge] (-1.5,-1) to (-1.5,1);
	\node[dot] at (0,0) {};
\end{scope}
\begin{scope}[shift={(-0.5,0)}]
	\draw[edge, out=90, in=60] (-0.5, -1) to (0,0);
	\draw[edge] (0, 0) to (0.5,1);
	\draw[edge] (0.5, -1) to (0,0);
	\draw[edge, out=-60, in=-90] (0, 0) to (-0.5,1);
	\draw[edge] (1.5,-1) to (1.5,1);
	\node[dot] at (0,0) {};
\end{scope}
\begin{scope}[shift={(0.5,-2)}]
	\draw[edge, out=90, in=60] (-0.5, -1) to (0,0);
	\draw[edge, out=60, in=-90] (0, 0) to (0.5,1);
	\draw[edge, out=90, in=-60] (0.5, -1) to (0,0);
	\draw[edge] (0, 0) to (-0.5,1);
	\draw[edge] (-1.5,-1) to (-1.5,1);
	\node[dot] at (0,0) {};
\end{scope}
\end{pgfonlayer}
\end{scope}
\end{tikzpicture}
\end{equation}
Noting that all cells of the same dimension have the same shape in $I^3$, we can eliminate all colouring by quotienting by the equivalence relation $E := \{(\sigma, \sigma') \;|\; d(\sigma) = d(\sigma')\}$. Models of $I^3/E$, that is, functors from the 3-category it generates to a monoidal category $\mathcal{C}$ (as a 2-category with a single 1-cell), are $R$-matrices in $\mathcal{C}$ \cite{perk2006yang}.

There is much more to quotients of cubes: in fact, they cover all sorts of \emph{associativity}-like equations. We introduce a couple general constructions for later use.
\begin{dfn}
Let $X$ be a computad. The \emph{cylinder} of X is the computad $I \otimes X$. The \emph{future cone} of $X$ is the quotient computad $C^+(X) := (I \otimes X)/(\{1\} \otimes X)$. The \emph{past cone} of $X$ is the quotient computad $C^-(X) := (I \otimes X)/(\{0\} \otimes X)$.
\end{dfn}
One can obtain different variants of these constructions by reversing the directions of cells of any dimension. All directed variants collapse for the corresponding constructions of CW complexes, the usual cylinder and cone of a topological space.

We can see the $n$-cube as the $n$-th iteration of the cylinder construction on the terminal computad. The corresponding iterations of the future cone produce Street's oriented simplexes, or \emph{orientals} \cite{street1987algebra}. This is proven in detail in \cite{buckley2015orientals}.
\begin{dfn}
The \emph{$n$-oriental} is the computad $C^+\stackrel{n}{\ldots} C^+(1)$.
\end{dfn}
The 2-oriental is just a binary operation.
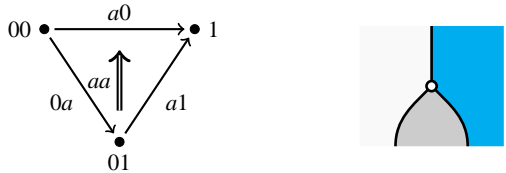
\begin{equation*}
\begin{tikzpicture}
\begin{pgfonlayer}{mid}
	\node[0c] (00) at (-1,0.5) [label=left:$00$] {}; 
	\node[0c] (01) at (0,-1) [label=below:$01$] {};
	\node[0c] (1) at (1,0.5) [label=right:$1$] {};
	\draw[1c] (00) to node[auto,swap] {$0a$} (01);
	\draw[1c] (01) to node[auto,swap] {$a1$} (1);
	\draw[1c] (00) to node[auto] {$a0$} (1);
	\draw[2c] (01) to node[auto] {$aa$} (0,0.5);
\end{pgfonlayer} 
\end{tikzpicture}
\hspace{50pt}
\begin{tikzpicture}[scale=0.8, baseline=-35pt]
\begin{pgfonlayer}{bg}
	\path[fill, color=gray!5] (-1.2,-1) -- (-0.6,-1) to [out=90, in=-135] (0,0) -- (0,1) -- (-1.2,1) -- cycle;
	\path[fill, color=gray!40] (-0.6,-1) to [out=90, in=-135] (0,0) to [out=-45, in= 90] (0.6,-1) -- cycle;
	\path[fill, color=cyan!100] (1.2,-1) -- (0.6,-1) to [out=90, in=-45] (0,0) to (0,1) -- (1.2,1) -- cycle;
\end{pgfonlayer}
\begin{pgfonlayer}{mid}
	\draw[edge, out=90, in=-135] (-0.6, -1) to (0,0);
	\draw[edge] (0, 0) to (0,1);
	\draw[edge, out=90, in=-45] (0.6, -1) to (0,0);
\end{pgfonlayer}
\begin{pgfonlayer}{fg}
	\node[dot] at (0,0) {};
\end{pgfonlayer}
\end{tikzpicture}
\end{equation*}

The 3-oriental is a coloured version of an associator.
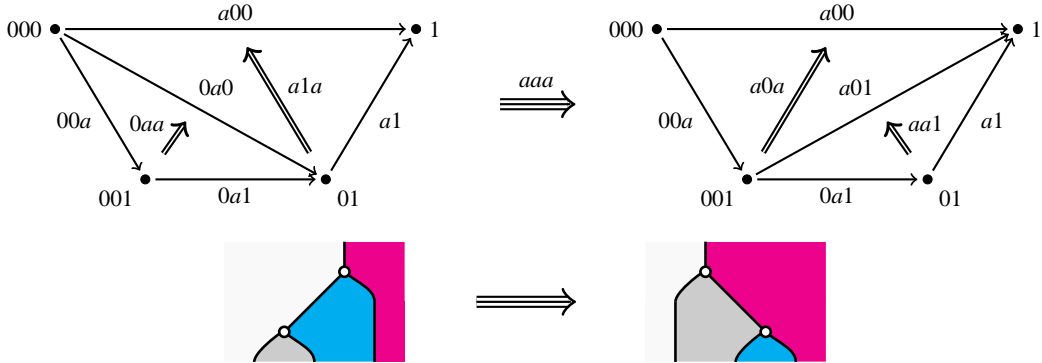
\begin{equation*}
\begin{tikzpicture}
\begin{pgfonlayer}{mid}
\begin{scope}[shift={(-4,0)}]
	\node[0c] (000) at (-2.4,1) [label=left:$000$] {}; 
	\node[0c] (001) at (-1.2,-1) [label=below left:$001$] {};
	\node[0c] (01) at (1.2,-1) [label=below right:$01$] {};
	\node[0c] (1) at (2.4,1) [label=right:$1$] {};
	\draw[1c] (000) to node[auto,swap] {$00a$} (001); 
	\draw[1c] (001) to node[auto,swap] {$0a1$} (01);
	\draw[1c] (01) to node[auto,swap] {$a1$} (1);
	\draw[1c] (000) to node[auto] {$a00$} (1);
	\draw[1c] (000) to node[auto] {$0a0$} (01); 
	\draw[2c] (001) to node[auto] {$0aa$} (-0.5,0); 
	\draw[2c] (01) to node[auto,swap] {$a1a$} (0,1);
\end{scope}
\begin{scope}
	\node at (0,0.3) {$aaa$};
	\draw[3c1] (-0.8,0) to (0.8,0);
	\draw[3c2] (-0.8,0) to (0.8,0);
	\draw[3c3] (-0.8,0) to (0.8,0);
\end{scope}
\begin{scope}[shift={(4,0)}]
	\node[0c] (000) at (-2.4,1) [label=left:$000$] {}; 
	\node[0c] (001) at (-1.2,-1) [label=below left:$001$] {};
	\node[0c] (01) at (1.2,-1) [label=below right:$01$] {};
	\node[0c] (1) at (2.4,1) [label=right:$1$] {};
	\draw[1c] (000) to node[auto,swap] {$00a$} (001); 
	\draw[1c] (001) to node[auto,swap] {$0a1$} (01);
	\draw[1c] (01) to node[auto,swap] {$a1$} (1);
	\draw[1c] (000) to node[auto] {$a00$} (1);
	\draw[1c] (001) to node[auto] {$a01$} (1); 
	\draw[2c] (01) to node[auto,swap] {$aa1$} (0.5,0); 
	\draw[2c] (001) to node[auto] {$a0a$} (0,1);
\end{scope}
\end{pgfonlayer} 
\end{tikzpicture}
\end{equation*}
\begin{equation*}
\begin{tikzpicture}[string, yscale=0.4, xscale=0.8]
\begin{scope}[shift={(-3.5,0)}]
\begin{pgfonlayer}{bg}
\begin{scope}[shift={(0.5,0.99)}]
	\path[fill, color=gray!5] (-2,-1) -- (-0.5,-1) to (0,0) to (0,1) -- (-2,1) -- cycle;
	\path[fill, color=cyan!100] (-0.5,-1) to (0,0) to [out=-60, in= 90] (0.5,-1) -- cycle;
	\path[fill, color=magenta!100] (1,-1) -- (0.5,-1) to [out=90, in=-60] (0,0) to (0,1) -- (1,1) -- cycle;
\end{scope}
\begin{scope}[shift={(-0.5,-0.99)}]
	\path[fill, color=gray!5] (-1,-1) -- (-0.5,-1) to [out=90, in=60] (0,0) -- (0.5,1) -- (-1,1) -- cycle;
	\path[fill, color=gray!40] (-0.5,-1) to [out=90, in=60] (0,0) to [out=-60, in= 90] (0.5,-1) -- cycle;
	\path[fill, color=cyan!100] (0.5,1) to (0,0) to [out=-60, in=90] (0.5,-1) -- (1.5,-1) -- (1.5,1) -- cycle;
	\path[fill, color=magenta!100] (1.5,-1) rectangle (2,1);
\end{scope}
\end{pgfonlayer}
\begin{pgfonlayer}{mid}
\begin{scope}[shift={(0.5,1)}]
	\draw[edge] (-0.5, -1) to (0,0);
	\draw[edge] (0, 0) to (0,1);
	\draw[edge, out=90, in=-60] (0.5, -1) to (0,0);
	\node[dot] at (0,0) {};
\end{scope}
\begin{scope}[shift={(-0.5,-1)}]
	\draw[edge, out=90, in=60] (-0.5, -1) to (0,0);
	\draw[edge] (0, 0) to (0.5,1);
	\draw[edge, out=90, in=-60] (0.5, -1) to (0,0);
	\draw[edge] (1.5,-1) to (1.5,1);
	\node[dot] at (0,0) {};
\end{scope}
\end{pgfonlayer}
\end{scope}
\begin{scope}
	\draw[3c1] (-1.2,0) to (1.2,0);
	\draw[3c2] (-1.2,0) to (1.2,0);
	\draw[3c3] (-1.2,0) to (1.2,0);
\end{scope}
\begin{scope}[shift={(3.5,0)}]
\begin{pgfonlayer}{bg}
\begin{scope}[shift={(-0.5,0.99)}]
	\path[fill, color=gray!5] (-1,-1) -- (-0.5,-1) to [out=90, in=60] (0,0) to (0,1) -- (-1,1) -- cycle;
	\path[fill, color=gray!40] (-0.5,-1) to [out=90, in=60] (0,0) to (0.5,-1) -- cycle;
	\path[fill, color=magenta!100] (2,-1) -- (0.5,-1) to (0,0) to (0,1) -- (2,1) -- cycle;
\end{scope}
\begin{scope}[shift={(0.5,-0.99)}]
	\path[fill, color=gray!5] (-2,-1) rectangle (-1.5,1);
	\path[fill, color=gray!40] (-1.5,-1) -- (-0.5,-1) to [out=90, in=60] (0,0) to (-0.5,1) -- (-1.5,1) -- cycle;
	\path[fill, color=cyan!100] (-0.5,-1) to [out=90, in=60] (0,0) to [out=-60, in= 90] (0.5,-1) -- cycle;
	\path[fill, color=magenta!100] (1,-1) -- (0.5,-1) to [out=90, in=-60] (0,0) to (-0.5,1) -- (1,1) -- cycle;
\end{scope}
\end{pgfonlayer}
\begin{pgfonlayer}{mid}
\begin{scope}[shift={(-0.5,1)}]
	\draw[edge, out=90, in=60] (-0.5, -1) to (0,0);
	\draw[edge] (0, 0) to (0,1);
	\draw[edge] (0.5, -1) to (0,0);
	\node[dot] at (0,0) {};
\end{scope}
\begin{scope}[shift={(0.5,-1)}]
	\draw[edge, out=90, in=60] (-0.5, -1) to (0,0);
	\draw[edge, out=90, in=-60] (0.5, -1) to (0,0);
	\draw[edge] (0, 0) to (-0.5,1);
	\draw[edge] (-1.5,-1) to (-1.5,1);
	\node[dot] at (0,0) {};
\end{scope}
\end{pgfonlayer}
\end{scope}
\end{tikzpicture}
\end{equation*}

The 4-oriental corresponds to the pentagonator, a directed version of MacLane's pentagon \cite{maclane1963natural} and a part of the theory of pseudomonoids; we only draw its string-diagrammatic version.
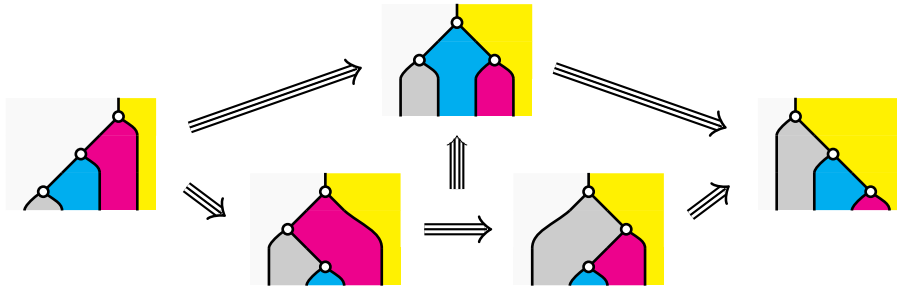
\begin{equation*}
\input{diagrams/pentagon.tikz}
\end{equation*}
Taking \emph{past} instead of future cones gives the coalgebraic duals of these constructions: respectively, a co-multiplication, co-associator, and co-pentagonator. 

Clearly, these computads can also be obtained directly as quotients of the $n$-cubes. Now, let $E$ instead be the equivalence relation on $I^3$ with
\begin{align*}
	& 0 \otimes \sigma \otimes \tau \sim 0 \otimes \sigma' \otimes \tau\;, \\
	& 1 \otimes \sigma \otimes \tau \sim 1 \otimes \sigma \otimes \tau'\;,
\end{align*}
for all $\sigma, \sigma', \tau, \tau' \in |I|$. Then, the quotient $I^3/E$ --- which roughly corresponds to taking a past cone on one copy of $I$, and a future cone on another, ``fibrewise'' on cells of the other copy --- presents a version of the \emph{Frobenius law} (see \cite[Chapter 4]{heunen2012lectures} for a review). This can be checked graphically by merging regions in Diagram (\ref{eq:yangbaxter}).
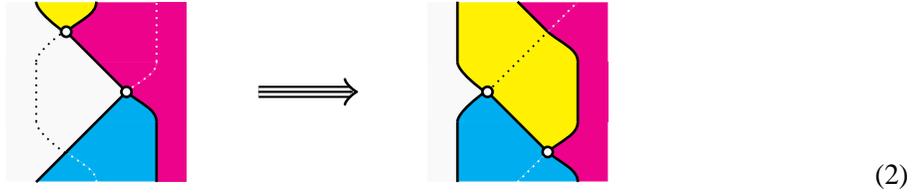
\begin{equation} \label{eq:frobenius}
\begin{tikzpicture}[string, yscale=0.4, xscale=0.8]
\begin{scope}[shift={(-3.5,0)}]
\begin{pgfonlayer}{bg}
\begin{scope}[shift={(-0.5,1.98)}]
	\path[fill, color=gray!5] (-1,-1) -- (-0.5,-1) to [out=90, in=60] (0,0) to [out=-60, in= -90] (-0.5,1) -- (-1,1) -- cycle;
	\path[fill, color=gray!5] (-0.5,-1) to [out=90, in=60] (0,0) to (0.5,-1) -- cycle;
	\path[fill, color=yellow!100] (-0.5,1) to [out=-90, in=-60] (0,0) to [out=60, in= -90] (0.5,1) -- cycle;
	\path[fill, color=magenta!100] (1.5,-1) -- (0.5,-1) to (0,0) to [out=60, in= -90] (0.5,1) -- (1.5,1) -- cycle;
	\path[fill, color=magenta!100] (1.5,-1) rectangle (2,1);
\end{scope}
\begin{scope}[shift={(0.5,0)}]
	\path[fill, color=gray!5] (-2,-1) rectangle (-1.5,1);
	\path[fill, color=gray!5] (-1.5,-1) -- (-0.5,-1) to (0,0) to (-0.5,1) -- (-1.5,1) -- cycle;
	\path[fill, color=cyan!100] (-0.5,-1) to (0,0) to [out=-60, in= 90] (0.5,-1) -- cycle;
	\path[fill, color=magenta!100] (-0.5,1) to (0,0) to [out=60, in= -90] (0.5,1) -- cycle;
	\path[fill, color=magenta!100] (1,-1) -- (0.5,-1) to [out=90, in=-60] (0,0) to [out=60, in= -90] (0.5,1) -- (1,1) -- cycle;
\end{scope}
\begin{scope}[shift={(-0.5,-1.98)}]
	\path[fill, color=gray!5] (-1,-1) -- (-0.5,-1) to (0,0) to [out=-60, in= -90] (-0.5,1) -- (-1,1) -- cycle;
	\path[fill, color=cyan!100] (-0.5,-1) to (0,0) to [out=-60, in= 90] (0.5,-1) -- cycle;
	\path[fill, color=gray!5] (-0.5,1) to [out=-90, in=-60] (0,0) to (0.5,1) -- cycle;
	\path[fill, color=cyan!100] (1.5,-1) -- (0.5,-1) to [out=90, in=-60] (0,0) to (0.5,1) -- (1.5,1) -- cycle;
	\path[fill, color=magenta!100] (1.5,-1) rectangle (2,1);
\end{scope}
\end{pgfonlayer}
\begin{pgfonlayer}{mid}
\begin{scope}[shift={(-0.5,2)}]
	\draw[edgedotdark, out=90, in=60] (-0.5, -1) to (0,0);
	\draw[edge, out=60, in=-90] (0, 0) to (0.5,1);
	\draw[edge] (0.5, -1) to (0,0);
	\draw[edge, out=-60, in=-90] (0, 0) to (-0.5,1);
	\draw[edgedot] (1.5,-1) to (1.5,1);
	\node[dot] at (0,0) {};
\end{scope}
\begin{scope}[shift={(0.5,0)}]
	\draw[edge] (-0.5, -1) to (0,0);
	\draw[edgedot, out=60, in=-90] (0, 0) to (0.5,1);
	\draw[edge, out=90, in=-60] (0.5, -1) to (0,0);
	\draw[edge] (0, 0) to (-0.5,1);
	\draw[edgedotdark] (-1.5,-1) to (-1.5,1);
	\node[dot] at (0,0) {};
\end{scope}
\begin{scope}[shift={(-0.5,-2)}]
	\draw[edgedot, out=90, in=-60] (0.5, -1) to (0,0);
	\draw[edgedotdark, out=-60, in=-90] (0, 0) to (-0.5,1);
	\draw[edge] (-0.5, -1) to (0,0);
	\draw[edge] (0, 0) to (0.5,1);
	\draw[edge] (1.5,-1) to (1.5,1);
\end{scope}
\end{pgfonlayer}
\end{scope}
\begin{scope}
	\draw[3c1] (-1.2,0) to (1.2,0);
	\draw[3c2] (-1.2,0) to (1.2,0);
	\draw[3c3] (-1.2,0) to (1.2,0);
\end{scope}
\begin{scope}[shift={(3.5,0)}]
\begin{pgfonlayer}{bg}
\begin{scope}[shift={(0.5,1.98)}]
	\path[fill, color=gray!5] (-2,-1) rectangle (-1.5,1);
	\path[fill, color=yellow!100] (-1.5,-1) -- (-0.5,-1) to (0,0) to (-0.5,1) -- (-1.5,1) -- cycle;
	\path[fill, color=yellow!100] (-0.5,-1) to (0,0) to [out=-60, in= 90] (0.5,-1) -- cycle;
	\path[fill, color=magenta!100] (-0.5,1) to (0,0) to (0.5,1) -- cycle;
	\path[fill, color=magenta!100] (1,-1) -- (0.5,-1) to [out=90, in=-60] (0,0) to (0.5,1) -- (1,1) -- cycle;
\end{scope}
\begin{scope}[shift={(-0.5,0)}]
	\path[fill, color=gray!5] (-1,-1) -- (-0.5,-1) to [out=90, in=60] (0,0) to [out=-60, in= -90] (-0.5,1) -- (-1,1) -- cycle;
	\path[fill, color=cyan!100] (-0.5,-1) to [out=90, in=60] (0,0) to (0.5,-1) -- cycle;
	\path[fill, color=yellow!100] (-0.5,1) to [out=-90, in=-60] (0,0) to (0.5,1) -- cycle;
	\path[fill, color=yellow!100] (1.5,-1) -- (0.5,-1) to (0,0) to (0.5,1) -- (1.5,1) -- cycle;
	\path[fill, color=magenta!100] (1.5,-1) rectangle (2,1);
\end{scope}
\begin{scope}[shift={(0.5,-1.98)}]
	\path[fill, color=gray!5] (-2,-1) rectangle (-1.5,1);
	\path[fill, color=cyan!100] (-1.5,-1) -- (-0.5,-1) to (0,0) to (-0.5,1) -- (-1.5,1) -- cycle;
	\path[fill, color=cyan!100] (-0.5,-1) to (0,0) to [out=-60, in= 90] (0.5,-1) -- cycle;
	\path[fill, color=yellow!100] (-0.5,1) to (0,0) to [out=60, in= -90] (0.5,1) -- cycle;
	\path[fill, color=magenta!100] (1,-1) -- (0.5,-1) to [out=90, in=-60] (0,0) to [out=60, in= -90] (0.5,1) -- (1,1) -- cycle;
\end{scope}
\end{pgfonlayer}
\begin{pgfonlayer}{mid}
\begin{scope}[shift={(0.5,2)}]
	\draw[edgedotdark] (-0.5, -1) to (0,0);
	\draw[edgedot] (0, 0) to (0.5,1);
	\draw[edge, out=90, in=-60] (0.5, -1) to (0,0);
	\draw[edge] (0, 0) to (-0.5,1);
	\draw[edge] (-1.5,-1) to (-1.5,1);
\end{scope}
\begin{scope}[shift={(-0.5,0)}]
	\draw[edge, out=90, in=60] (-0.5, -1) to (0,0);
	\draw[edgedotdark] (0, 0) to (0.5,1);
	\draw[edge] (0.5, -1) to (0,0);
	\draw[edge, out=-60, in=-90] (0, 0) to (-0.5,1);
	\draw[edge] (1.5,-1) to (1.5,1);
	\node[dot] at (0,0) {};
\end{scope}
\begin{scope}[shift={(0.5,-2)}]
	\draw[edgedot] (-0.5, -1) to (0,0);
	\draw[edge, out=60, in=-90] (0, 0) to (0.5,1);
	\draw[edge, out=90, in=-60] (0.5, -1) to (0,0);
	\draw[edge] (0, 0) to (-0.5,1);
	\draw[edge] (-1.5,-1) to (-1.5,1);
	\node[dot] at (0,0) {};
\end{scope}
\end{pgfonlayer}
\end{scope}
\end{tikzpicture}
\end{equation}

Again, since both orientals and their duals have congruent shapes for cells of equal dimension, we can quotient in order to obtain presentations of the theory of semigroups and co-semigroups, as lax --- that is, with as many non-invertible higher cells --- as we want them to be.

Now, we show how to obtain the theory of \emph{monoids} by cones, in two steps; for comonoids, it will suffice to dualise everything. First of all, we need the 2-computad $K$ presenting the ``theory of constants''.
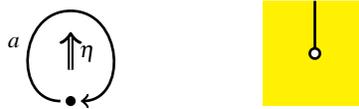
\begin{equation*}
\begin{tikzpicture}
\begin{pgfonlayer}{mid}
	\node[0c] (0) at (0,-1.2) {};
	\draw[edge, shorten <=2pt] (0) to [out=180, in=180, looseness=1.5] node[auto] {$a$} (0,0);
	\draw[edge, ->, shorten >=2pt] (0,0) to [out=0, in=0, looseness=1.5] (0);
	\draw[2c] (0) to node[auto,swap] {$\eta$} (0,0);
\end{pgfonlayer} 
\end{tikzpicture}
\hspace{50pt}
\begin{tikzpicture}[string, scale=0.7]
\begin{pgfonlayer}{bg}
	\fill[color=yellow!100] (-1,-1) rectangle (1,1);
\end{pgfonlayer}
\begin{pgfonlayer}{mid}
	\draw[edge] (0,0) to (0,1);
	\node[dot] at (0,0) {};
\end{pgfonlayer}
\end{tikzpicture}
\end{equation*}
Then, the future cone of $K$ with reversed $1$-cells, $C^+(K)^\mathrm{op(1)}$, is a 3-computad presenting an operation with a right unit.
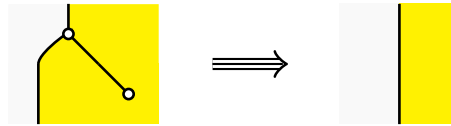
\begin{equation*}
\begin{tikzpicture}[string, yscale=0.4, xscale=0.8]
\begin{scope}[shift={(-2.5,0)}] 
\begin{pgfonlayer}{bg}
\begin{scope}[shift={(-0.5,0.99)}]
	\path[fill, color=gray!5] (-1,-1) -- (-0.5,-1) to [out=90, in=60] (0,0) to (0,1) -- (-1,1) -- cycle;
	\path[fill, color=yellow!100] (-0.5,-1) to [out=90, in=60] (0,0) -- (0,1) -- (1.5,1) -- (1.5,-1) -- cycle;
\end{scope}
\begin{scope}[shift={(0.5,-0.99)}]
	\path[fill, color=gray!5] (-2,-1) rectangle (-1.5,1);
	\path[fill, color=yellow!100] (-1.5,-1) rectangle (0.5,1);
\end{scope}
\end{pgfonlayer}
\begin{pgfonlayer}{mid}
\begin{scope}[shift={(-0.5,1)}]
	\draw[edge, out=90, in=60] (-0.5, -1) to (0,0);
	\draw[edge] (0, 0) to (0,1);
	\draw[edge] (0.5, -1) to (0,0);
	\node[dot] at (0,0) {};
\end{scope}
\begin{scope}[shift={(0.5,-1)}]
	\draw[edge] (0, 0) to (-0.5,1);
	\draw[edge] (-1.5,-1) to (-1.5,1);
	\node[dot] at (0,0) {};
\end{scope}
\end{pgfonlayer}
\end{scope}
\begin{scope}[shift={(2.5,0)}] 
\begin{pgfonlayer}{bg}
\begin{scope}
	\path[fill, color=gray!5] (-1,-2) rectangle (0,2);
	\path[fill, color=yellow!100] (0,-2) rectangle (1,2);
\end{scope}
\end{pgfonlayer}
\begin{pgfonlayer}{mid}
\begin{scope}
	\draw[edge] (0, -2) to (0,2);
\end{scope}
\end{pgfonlayer}
\end{scope}
\begin{scope}
	\draw[3c1] (-1,0) to (1,0);
	\draw[3c2] (-1,0) to (1,0);
	\draw[3c3] (-1,0) to (1,0);
\end{scope}
\end{tikzpicture}
\end{equation*}
Now, take the future cone $C^+(C^+(K)^\mathrm{op(1)})$; this can be pictured as follows.
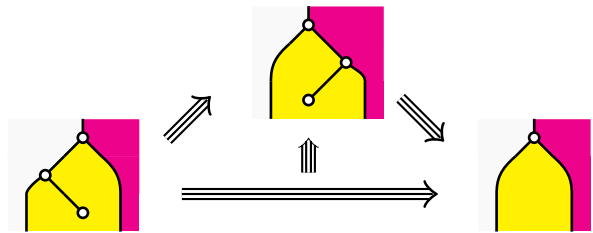
\begin{equation*}
\begin{tikzpicture}[string, yscale=0.25, xscale=0.5]

\begin{scope}[shift={(-6,-5)}] 
\begin{pgfonlayer}{bg}
\begin{scope}[shift={(0,1.98)}]
	\path[fill, color=gray!5] (-2,-1) -- (-0.5,-1) to (0,0) to (0,1) -- (-2,1) -- cycle;
	\path[fill, color=yellow!100] (-0.5,-1) to (0,0) to (0.5,-1) -- cycle;
	\path[fill, color=magenta!100] (1.5,-1) -- (0.5,-1) to (0,0) to (0,1) -- (1.5,1) -- cycle;
\end{scope}
\begin{scope}[shift={(-1,0)}]
	\path[fill, color=gray!5] (-1,-1) -- (-0.5,-1) to [out=90, in=60] (0,0) to (0.5,1) -- (-1,1) -- cycle;
	\path[fill, color=yellow!100] (-0.5,-1) to [out=90, in=60] (0,0) to (0.5,1) -- (1.5,1) to [out=-60, in=90] (2,-1) -- cycle;
	\path[fill, color=magenta!100] (1.5,1) to [out=-60, in=90] (2,-1) -- (2.5,-1) -- (2.5,1) -- cycle;
\end{scope}
\begin{scope}[shift={(0,-1.98)}]
	\path[fill, color=gray!5] (-2,-1) rectangle (-1.5,1);
	\path[fill, color=yellow!100] (-1.5,-1) rectangle (1,1);
	\path[fill, color=magenta!100] (1,-1) rectangle (1.5,1);
\end{scope}
\end{pgfonlayer}
\begin{pgfonlayer}{mid}
\begin{scope}[shift={(0,2)}]
	\draw[edge] (-0.5, -1) to (0,0);
	\draw[edge] (0, 0) to (0,1);
	\draw[edge] (0.5, -1) to (0,0);
	\node[dot] at (0,0) {};
\end{scope}
\begin{scope}[shift={(-1,0)}]
	\draw[edge, out=90, in=60] (-0.5, -1) to (0,0);
	\draw[edge] (0, 0) to (0.5,1);
	\draw[edge] (0.5, -1) to (0,0);
	\draw[edge, out=-60, in=90] (1.5,1) to (2,-1);
	\node[dot] at (0,0) {};
\end{scope}
\begin{scope}[shift={(0,-2)}]
	\draw[edge] (0, 0) to (-0.5,1);
	\draw[edge] (-1.5,-1) to (-1.5,1);
	\draw[edge] (1,-1) to (1,1);
	\node[dot] at (0,0) {};
\end{scope}
\end{pgfonlayer}
\end{scope}
\begin{scope}[shift={(0,1)}] 
\begin{pgfonlayer}{bg}
\begin{scope}[shift={(0,1.98)}]
	\path[fill, color=gray!5] (-1.5,-1) -- (-0.5,-1) to (0,0) to (0.5,1) -- (-1.5,1) -- cycle;
	\path[fill, color=yellow!100] (-0.5,-1) to (0,0) to (0.5,-1) -- cycle;
	\path[fill, color=magenta!100] (2,-1) -- (0.5,-1) to (0,0) to (0,1) -- (2,1) -- cycle;
\end{scope}
\begin{scope}[shift={(1,0)}]
	\path[fill, color=gray!5] (-2.5,-1) -- (-2,-1) to [out=90, in=60] (-1.5,1) -- (-2.5,1) -- cycle;
	\path[fill, color=yellow!100] (-2,-1) to [out=90, in=-120] (-1.5,1) -- (-0.5,1) -- (0,0) to [out=-60, in=90] (0.5,-1) -- (-2,-1) -- cycle;
	\path[fill, color=magenta!100] (1,-1) -- (0.5,-1) to [out=90, in=-60] (0,0) to (-0.5,1) -- (1,1) -- cycle;
\end{scope}
\begin{scope}[shift={(0,-1.98)}]
	\path[fill, color=gray!5] (-1.5,-1) rectangle (-1,1);
	\path[fill, color=yellow!100] (-1,-1) rectangle (1.5,1);
	\path[fill, color=magenta!100] (1.5,-1) rectangle (2,1);
\end{scope}
\end{pgfonlayer}
\begin{pgfonlayer}{mid}
\begin{scope}[shift={(0,2)}]
	\draw[edge] (-0.5, -1) to (0,0);
	\draw[edge] (0, 0) to (0,1);
	\draw[edge] (0.5, -1) to (0,0);
	\node[dot] at (0,0) {};
\end{scope}
\begin{scope}[shift={(1,0)}]
	\draw[edge] (-0.5, -1) to (0,0);
	\draw[edge, out=90, in=-60] (0.5, -1) to (0,0);
	\draw[edge] (0, 0) to (-0.5,1);
	\draw[edge, out=90, in=-120] (-2,-1) to (-1.5,1);
	\node[dot] at (0,0) {};
\end{scope}
\begin{scope}[shift={(0,-2)}]
	\draw[edge] (0, 0) to (0.5,1);
	\draw[edge] (-1,-1) to (-1,1);
	\draw[edge] (1.5,-1) to (1.5,1);
	\node[dot] at (0,0) {};
\end{scope}
\end{pgfonlayer}
\end{scope}

\begin{scope}[shift={(6,-5)}] 
\begin{pgfonlayer}{bg}
\begin{scope}[shift={(0,1.98)}]
	\path[fill, color=gray!5] (-1.5,-1) -- (-0.5,-1) to (0,0) to (0.5,1) -- (-1.5,1) -- cycle;
	\path[fill, color=yellow!100] (-0.5,-1) to (0,0) to (0.5,-1) -- cycle;
	\path[fill, color=magenta!100] (1.5,-1) -- (0.5,-1) to (0,0) to (0,1) -- (1.5,1) -- cycle;
\end{scope}
\begin{scope}[shift={(1,0)}]
	\path[fill, color=gray!5] (-2.5,-1) rectangle (-1.5,1);
	\path[fill, color=magenta!100] (-0.5,-1) rectangle (0.5,1);
	\path[fill, color=yellow!100] (-2,-1) to [out=90, in=-120] (-1.5,1) -- (-0.5,1) to [out=-60, in=90] (0,-1) -- cycle;
\end{scope}
\begin{scope}[shift={(0,-1.98)}]
	\path[fill, color=gray!5] (-1.5,-1) rectangle (-1,1);
	\path[fill, color=yellow!100] (-1,-1) rectangle (1,1);
	\path[fill, color=magenta!100] (1,-1) rectangle (1.5,1);
\end{scope}
\end{pgfonlayer}
\begin{pgfonlayer}{mid}
\begin{scope}[shift={(0,2)}]
	\draw[edge] (-0.5, -1) to (0,0);
	\draw[edge] (0, 0) to (0,1);
	\draw[edge] (0.5, -1) to (0,0);
	\node[dot] at (0,0) {};
\end{scope}
\begin{scope}[shift={(1,0)}]
	\draw[edge, out=90, in=-120] (-2, -1) to (-1.5,1);
	\draw[edge, out=90, in=-60] (0, -1) to (-0.5,1);
\end{scope}
\begin{scope}[shift={(0,-2)}]
	\draw[edge] (-1,-1) to (-1,1);
	\draw[edge] (1,-1) to (1,1);
\end{scope}
\end{pgfonlayer}
\end{scope}
\begin{scope}[shift={(0,-6)}] 
	\draw[3c1] (-4,0) to (4,0);
	\draw[3c2] (-4,0) to (4,0);
	\draw[3c3] (-4,0) to (4,0);
\end{scope}
\begin{scope}[shift={(3,0)}] 
	\draw[3c1] (-1,0) to (1,-4);
	\draw[3c2] (-1,0) to (1,-4);
	\draw[3c3] (-1,0) to (1,-4);
\end{scope}
\begin{scope}[shift={(-3.2,0)}] 
	\draw[3c1] (-1,-4) to (1,0);
	\draw[3c2] (-1,-4) to (1,0);
	\draw[3c3] (-1,-4) to (1,0);
\end{scope}
\begin{scope}[shift={(0,-3)}] 
	\draw[4c1] (0,-2) to (0,0);
	\draw[4c2] (0,-2) to (0,0);
\end{scope}
\end{tikzpicture}
\end{equation*}
This is a directed version of MacLane's triangle equality. As usual, we can identify all cells of congruent shape; since we have a right unitor 3-cell for the ``white-yellow'' operation, and a new left unitor 3-cell for the ``yellow-magenta'' operation, upon identification, the constant becomes a two-sided unit. 

Thus, by quotienting $C^+(C^+(K)^\mathrm{op(1)})$ (or further iterations of future cones), we obtain presentations of arbitrarily lax theories of monoids.

\section{Homomorphisms and actions} \label{sec:homo}

Now that we have a list of basic examples, let us reconsider the constructions of Section \ref{sec:basic} more in general. Given a computad $X$, the cylinder of $X$ has the following structure:
\begin{enumerate}
	\item it contains two copies of $X$, $\{0\} \otimes X$ and $\{1\} \otimes X$;
	\item for all generating $n$-cells $\sigma$ of $X$, it contains an $(n+1)$-cell $a \otimes \sigma$, which has $0 \otimes \sigma$ in its input border, and $1 \otimes \sigma$ in its output border.
\end{enumerate}
The category $\omega\mathbf{Cat}$ with the Crans-Gray tensor product actually admits a \emph{biclosed} structure, so that a functor $h: \mathcal{I} \otimes \mathcal{X} \to \mathcal{C}$ corresponds to a functor $\tilde{h}: \mathcal{I} \to (\mathcal{X} \multimap \mathcal{C})$, where $\mathcal{X} \multimap \mathcal{C}$ is an internal hom-object. But $\mathcal{I}$ is just the ``walking arrow'' category, so $\tilde{h}$ is a morphism in $\mathcal{X} \multimap \mathcal{C}$ --- which can be seen as a higher-dimensional version of a category of $X$-algebras in $\mathcal{C}$; hence, $\tilde{h}$ is a \emph{homomorphism} of $X$-algebras.

This suggests that $I \otimes X$ captures the theory of homomorphisms for the theory $X$: a functor $h: \mathcal{I} \otimes \mathcal{X} \to \mathcal{C}$ is an internal homomorphism in $\mathcal{C}$ between the two models $h_0: \{0\} \otimes \mathcal{X} \to \mathcal{C}$ and $h_1: \{1\} \otimes \mathcal{X} \to \mathcal{C}$. Observe that the undirected analogue is just a \emph{homotopy} of maps of topological spaces.

\begin{exm}
Let $M$ be a computad presenting a theory of (possibly lax) monoids, with a single ``object'' 1-cell $a$, a multiplication 2-cell $\mu$, and a unit 2-cell $\eta$. In $I \otimes M$, the 3-cells $a \otimes \mu$ and $a \otimes \eta$ can be pictured in the following way.
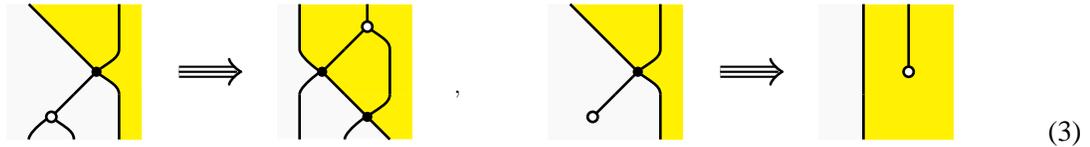
\begin{equation} \label{eq:monoids}
\begin{tikzpicture}[string, yscale=0.3, xscale=0.6]
\begin{scope}[shift={(-3,0)}]
\begin{pgfonlayer}{bg}
\begin{scope}[shift={(-0.5,1.99)}]
	\path[fill, color=gray!5] (-1,-1) -- (0.5,-1) -- (-0.5, 1) -- (-1,1) -- cycle;
	\path[fill, color=yellow!100] (-0.5,1) -- (0.5, -1) -- (2,-1) -- (2,1) -- cycle;
\end{scope}
\begin{scope}[shift={(0.5,0)}]
	\path[fill, color=gray!5] (-2,-1) -- (0.5,-1) to [out=90, in=-60] (0,0) -- (-0.5,1) -- (-2,1) -- cycle;
	\path[fill, color=yellow!100] (-0.5,1) to (0,0) to [out=-60, in=90] (0.5,-1) -- (1,-1) -- (1,1) --  cycle;
\end{scope}
\begin{scope}[shift={(-0.5,-1.99)}]
	\path[fill, color=gray!5] (-1,-1) rectangle (1.5,1);
	\path[fill, color=yellow!100] (1.5,-1) rectangle (2,1);
\end{scope}
\end{pgfonlayer}
\begin{pgfonlayer}{mid}
\begin{scope}[shift={(-0.5,2)}]
	\draw[edge] (0.5, -1) to (-0.5,1);
	\draw[edge] (1.5,-1) to (1.5,1);
\end{scope}
\begin{scope}[shift={(0.5,0)}]
	\draw[edge] (-0.5, -1) to (0,0);
	\draw[edge] (-0.5, 1) to (0,0);
	\draw[edge, out=60, in=-90] (0, 0) to (0.5,1);
	\draw[edge, out=90, in=-60] (0.5, -1) to (0,0);
	\node[dotdark] at (0,0) {};
\end{scope}
\begin{scope}[shift={(-0.5,-2)}]
	\draw[edge, out=90, in=60] (-0.5, -1) to (0,0);
	\draw[edge] (0, 0) to (0.5,1);
	\draw[edge, out=90, in=-60] (0.5, -1) to (0,0);
	\draw[edge] (1.5,-1) to (1.5,1);
	\node[dot] at (0,0) {};
\end{scope}
\end{pgfonlayer}
\end{scope}
\begin{scope}
	\draw[3c1] (-1.2,0) to (1.2,0);
	\draw[3c2] (-1.2,0) to (1.2,0);
	\draw[3c3] (-1.2,0) to (1.2,0);
\end{scope}
\begin{scope}[shift={(3,0)}]
\begin{pgfonlayer}{bg}
\begin{scope}[shift={(0.5,1.99)}]
	\path[fill, color=gray!5] (-2,-1) rectangle (-1.5,1);
	\path[fill, color=yellow!100] (-1.5,-1) rectangle (1,1);
\end{scope}
\begin{scope}[shift={(-0.5,0)}]
	\path[fill, color=gray!5] (-1,-1) -- (0.5,-1) to (0,0) to [out=-60, in= -90] (-0.5,1) -- (-1,1) -- cycle;
	\path[fill, color=yellow!100] (-0.5,1) to [out=-90, in=-60] (0,0) to (0.5,-1) -- (2,-1) -- (2,1) -- cycle;
\end{scope}
\begin{scope}[shift={(0.5,-1.99)}]
	\path[fill, color=gray!5] (-2,-1) -- (0.5,-1) -- (-0.5,1) -- (-2,1) -- cycle;
	\path[fill, color=yellow!100] (-0.5,1) -- (1,1) -- (1,-1) -- (0.5,-1) -- cycle;
\end{scope}
\end{pgfonlayer}
\begin{pgfonlayer}{mid}
\begin{scope}[shift={(0.5,2)}]
	\draw[edge] (-0.5, -1) to (0,0);
	\draw[edge] (0, 0) to (0,1);
	\draw[edge, out=90, in=-60] (0.5, -1) to (0,0);
	\draw[edge] (-1.5,-1) to (-1.5,1);
	\node[dot] at (0,0) {};
\end{scope}
\begin{scope}[shift={(-0.5,0)}]
	\draw[edge, out=90, in=60] (-0.5, -1) to (0,0);
	\draw[edge] (0, 0) to (0.5,1);
	\draw[edge] (0.5, -1) to (0,0);
	\draw[edge, out=-60, in=-90] (0, 0) to (-0.5,1);
	\draw[edge] (1.5,-1) to (1.5,1);
	\node[dotdark] at (0,0) {};
\end{scope}
\begin{scope}[shift={(0.5,-2)}]
	\draw[edge, out=90, in=60] (-0.5, -1) to (0,0);
	\draw[edge, out=60, in=-90] (0, 0) to (0.5,1);
	\draw[edge] (0.5, -1) to (0,0);
	\draw[edge] (-1.5,-1) to (-1.5,1);
	\draw[edge] (0, 0) to (-0.5,1);
	\node[dotdark] at (0,0) {};
\end{scope}
\end{pgfonlayer}
\end{scope}
\begin{scope}[shift={(5.5,-1)}] 
	\node at (0,0) {$,$};
\end{scope}
\end{tikzpicture}
\hspace{30pt}
\begin{tikzpicture}[string, yscale=0.3, xscale=0.6]
\begin{scope}[shift={(-3,0)}]
\begin{pgfonlayer}{bg}
\begin{scope}[shift={(-0.5,1.99)}]
	\path[fill, color=gray!5] (-1,-1) -- (0.5,-1) -- (-0.5, 1) -- (-1,1) -- cycle;
	\path[fill, color=yellow!100] (-0.5,1) -- (0.5, -1) -- (2,-1) -- (2,1) -- cycle;
\end{scope}
\begin{scope}[shift={(0.5,0)}]
	\path[fill, color=gray!5] (-2,-1) -- (0.5,-1) to [out=90, in=-60] (0,0) -- (-0.5,1) -- (-2,1) -- cycle;
	\path[fill, color=yellow!100] (-0.5,1) to (0,0) to [out=-60, in=90] (0.5,-1) -- (1,-1) -- (1,1) --  cycle;
\end{scope}
\begin{scope}[shift={(-0.5,-1.99)}]
	\path[fill, color=gray!5] (-1,-1) rectangle (1.5,1);
	\path[fill, color=yellow!100] (1.5,-1) rectangle (2,1);
\end{scope}
\end{pgfonlayer}
\begin{pgfonlayer}{mid}
\begin{scope}[shift={(-0.5,2)}]
	\draw[edge] (0.5, -1) to (-0.5,1);
	\draw[edge] (1.5,-1) to (1.5,1);
\end{scope}
\begin{scope}[shift={(0.5,0)}]
	\draw[edge] (-0.5, -1) to (0,0);
	\draw[edge] (-0.5, 1) to (0,0);
	\draw[edge, out=60, in=-90] (0, 0) to (0.5,1);
	\draw[edge, out=90, in=-60] (0.5, -1) to (0,0);
	\node[dotdark] at (0,0) {};
\end{scope}
\begin{scope}[shift={(-0.5,-2)}]
	\draw[edge] (0, 0) to (0.5,1);
	\draw[edge] (1.5,-1) to (1.5,1);
	\node[dot] at (0,0) {};
\end{scope}
\end{pgfonlayer}
\end{scope}
\begin{scope}
	\draw[3c1] (-1.2,0) to (1.2,0);
	\draw[3c2] (-1.2,0) to (1.2,0);
	\draw[3c3] (-1.2,0) to (1.2,0);
\end{scope}
\begin{scope}[shift={(3,0)}]
\begin{pgfonlayer}{bg}
\begin{scope}[shift={(0.5,1.99)}]
	\path[fill, color=gray!5] (-2,-1) rectangle (-1,1);
	\path[fill, color=yellow!100] (-1,-1) rectangle (1,1);
\end{scope}
\begin{scope}[shift={(0.5,0)}]
	\path[fill, color=gray!5] (-2,-1) rectangle (-1,1);
	\path[fill, color=yellow!100] (-1,-1) rectangle (1,1);
\end{scope}
\begin{scope}[shift={(0.5,-1.99)}]
	\path[fill, color=gray!5] (-2,-1) rectangle (-1,1);
	\path[fill, color=yellow!100] (-1,-1) rectangle (1,1);
\end{scope}
\end{pgfonlayer}
\begin{pgfonlayer}{mid}
\begin{scope}[shift={(0.5,2)}]
	\draw[edge] (-1,-1) to (-1,1);
	\draw[edge] (0,1) to (0,-1);
\end{scope}
\begin{scope}[shift={(0.5,0)}]
	\draw[edge] (-1,-1) to (-1,1);
	\draw[edge] (0,0) to (0,1);
	\node[dot] at (0,0) {};
\end{scope}
\begin{scope}[shift={(0.5,-2)}]
	\draw[edge] (-1,-1) to (-1,1);
\end{scope}
\end{pgfonlayer}
\end{scope}
\end{tikzpicture}
\end{equation}
These cells embody the action of ``sliding'' the multiplication and unit past the mediating 1-cell $a \otimes *$, where $*$ is the unique $0$-cell of $M$. The use of sliding for reasoning about naturality is discussed extensively in \cite{hinze2016equational}; the tensor product of computads provides it with a compositional semantics.

Models of $I \otimes M$ in the 2-category $\mathbf{Cat}$ are pairs of monads related by a Kleisli law, an asymmetric version of a \emph{distributive law} \cite{beck1969distributive}. The dual notion of Eilenberg-Moore law is obtained as a model of $(I \otimes M)^{\mathrm{op}(1)}$.
\end{exm}

If $X$ is PRO-like, in the sense that it has only one 0-cell $*$ and one cell $a$ of lowest, non-zero dimension $k$ (so that $k$ is seen as the dimension of ``objects'' in the theory), one is usually interested in homomorphisms between $X$-algebras in a higher category $\mathcal{C}$ that are localised at the same 0-cell, so that the lowest-dimensional component is a $(k+1)$-morphism between the two underlying objects of the algebras. 

In order to achieve this, we quotient $I \otimes X$ to obtain the \emph{reduced cylinder} of $X$, $(I \otimes X)/(I \otimes \{*\})$. Again, this construction has a well-established undirected analogue: maps from the reduced cylinder of a topological space into another space are pointed homotopies.

\begin{exm}
For the computad $M$ of the previous example, the 3-cells $a \otimes \mu$ and $a \otimes \eta$ take the following form in the reduced cylinder of $M$.
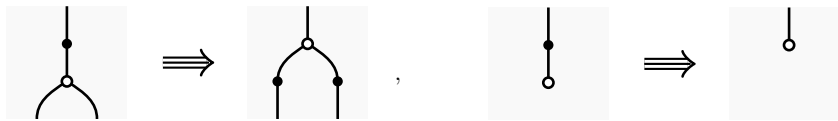
\begin{equation*}
\begin{tikzpicture}[xscale=0.8, yscale=0.5]
\begin{scope}[shift={(-2,0)}] 
\begin{pgfonlayer}{bg}
	\path[fill, color=gray!5] (-1,-2) rectangle (1,1);
\end{pgfonlayer}
\begin{pgfonlayer}{mid}
	\draw[edge, out=90, in=-135] (-0.5, -2) to (0,-1);
	\draw[edge] (0, -1) to (0,1);
	\draw[edge, out=90, in=-45] (0.5, -2) to (0,-1);
	\node[dot] at (0,-1) {};
	\node[dotdark] at (0,0) {};
\end{pgfonlayer}
\end{scope}
\begin{scope}[shift={(0,-0.5)}] 
	\draw[3c1] (-0.8,0) to (0.8,0);
	\draw[3c2] (-0.8,0) to (0.8,0);
	\draw[3c3] (-0.8,0) to (0.8,0);
\end{scope}
\begin{scope}[shift={(2,0)}] 
\begin{pgfonlayer}{bg}
	\path[fill, color=gray!5] (-1,-2) rectangle (1,1);
\end{pgfonlayer}
\begin{pgfonlayer}{mid}
	\draw[edge, out=90, in=-135] (-0.5, -1) to (0,0);
	\draw[edge] (0, 0) to (0,1);
	\draw[edge, out=90, in=-45] (0.5, -1) to (0,0);
	\draw[edge] (-0.5, -2) to (-0.5, -1);
	\draw[edge] (0.5, -2) to (0.5, -1);
	\node[dot] at (0,0) {};
	\node[dotdark] at (-0.5,-1) {};
	\node[dotdark] at (0.5,-1) {};
\end{pgfonlayer}
\end{scope}
\begin{scope}[shift={(3.5,-1)}] 
	\node at (0,0) {$,$};
\end{scope}
\end{tikzpicture}
\hspace{30pt}
\begin{tikzpicture}[xscale=0.8, yscale=0.5]
\begin{scope}[shift={(-2,0)}] 
\begin{pgfonlayer}{bg}
	\path[fill, color=gray!5] (-1,-2) rectangle (1,1);
\end{pgfonlayer}
\begin{pgfonlayer}{mid}
	\draw[edge] (0, -1) to (0,1);
	\node[dot] at (0,-1) {};
	\node[dotdark] at (0,0) {};
\end{pgfonlayer}
\end{scope}
\begin{scope}[shift={(0,-0.5)}] 
	\draw[3c1] (-0.8,0) to (0.8,0);
	\draw[3c2] (-0.8,0) to (0.8,0);
	\draw[3c3] (-0.8,0) to (0.8,0);
\end{scope}
\begin{scope}[shift={(2,0)}] 
\begin{pgfonlayer}{bg}
	\path[fill, color=gray!5] (-1,-2) rectangle (1,1);
\end{pgfonlayer}
\begin{pgfonlayer}{mid}
	\draw[edge] (0, 0) to (0,1);
	\node[dot] at (0,0) {};
\end{pgfonlayer}
\end{scope}
\end{tikzpicture}
\end{equation*}
\end{exm}

In the \emph{cone} of a computad $X$, one of the copies of $X$ in a cylinder is trivialised by a quotient. After this analysis, the future cone of $X$ can be seen as the theory of homomorphisms from an arbitrary $X$-algebra to the trivial $X$-algebra. When $X$ is, for instance, the theory of monoids, this happens to capture precisely the notion of \emph{left action} of a monoid, as can most easily be seen by merging the regions beyond $a \otimes *$ in Diagram (\ref{eq:monoids}).
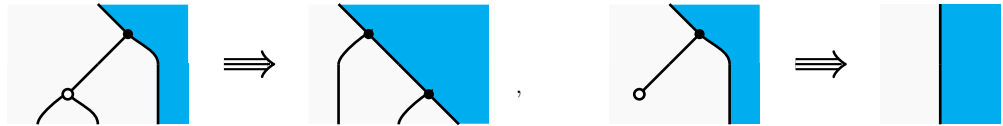
\begin{equation*}
\begin{tikzpicture}[string, yscale=0.4, xscale=0.8]
\begin{scope}[shift={(-2.5,0)}]
\begin{pgfonlayer}{bg}
\begin{scope}[shift={(0.5,0.98)}]
	\path[fill, color=gray!5] (-2,-1) -- (0.5,-1) to [out=90, in=-60] (0,0) -- (-0.5,1) -- (-2,1) -- cycle;
	\path[fill, color=cyan!100] (-0.5,1) to (0,0) to [out=-60, in=90] (0.5,-1) -- (1,-1) -- (1,1) --  cycle;
\end{scope}
\begin{scope}[shift={(-0.5,-0.98)}]
	\path[fill, color=gray!5] (-1,-1) rectangle (1.5,1);
	\path[fill, color=cyan!100] (1.5,-1) rectangle (2,1);
\end{scope}
\end{pgfonlayer}
\begin{pgfonlayer}{mid}
\begin{scope}[shift={(0.5,1)}]
	\draw[edge] (-0.5, -1) to (0,0);
	\draw[edge] (-0.5, 1) to (0,0);
	\draw[edge, out=90, in=-60] (0.5, -1) to (0,0);
	\node[dotdark] at (0,0) {};
\end{scope}
\begin{scope}[shift={(-0.5,-1)}]
	\draw[edge, out=90, in=60] (-0.5, -1) to (0,0);
	\draw[edge] (0, 0) to (0.5,1);
	\draw[edge, out=90, in=-60] (0.5, -1) to (0,0);
	\draw[edge] (1.5,-1) to (1.5,1);
	\node[dot] at (0,0) {};
\end{scope}
\end{pgfonlayer}
\end{scope}
\begin{scope}
	\draw[3c1] (-0.8,0) to (0.8,0);
	\draw[3c2] (-0.8,0) to (0.8,0);
	\draw[3c3] (-0.8,0) to (0.8,0);
\end{scope}
\begin{scope}[shift={(2.5,0)}]
\begin{pgfonlayer}{bg}
\begin{scope}[shift={(-0.5,0.98)}]
	\path[fill, color=gray!5] (-1,-1) -- (0.5,-1) to (0,0) to (-0.5,1) -- (-1,1) -- cycle;
	\path[fill, color=cyan!100] (-0.5,1) to (0,0) to (0.5,-1) -- (2,-1) -- (2,1) -- cycle;
\end{scope}
\begin{scope}[shift={(0.5,-0.98)}]
	\path[fill, color=gray!5] (-2,-1) -- (0.5,-1) -- (-0.5,1) -- (-2,1) -- cycle;
	\path[fill, color=cyan!100] (-0.5,1) -- (1,1) -- (1,-1) -- (0.5,-1) -- cycle;
\end{scope}
\end{pgfonlayer}
\begin{pgfonlayer}{mid}
\begin{scope}[shift={(-0.5,1)}]
	\draw[edge, out=90, in=60] (-0.5, -1) to (0,0);
	\draw[edge] (0.5, -1) to (0,0);
	\draw[edge] (0, 0) to (-0.5,1);
	\node[dotdark] at (0,0) {};
\end{scope}
\begin{scope}[shift={(0.5,-1)}]
	\draw[edge, out=90, in=60] (-0.5, -1) to (0,0);
	\draw[edge] (0.5, -1) to (0,0);
	\draw[edge] (-1.5,-1) to (-1.5,1);
	\draw[edge] (0, 0) to (-0.5,1);
	\node[dotdark] at (0,0) {};
\end{scope}
\end{pgfonlayer}
\end{scope}
\begin{scope}[shift={(4.5,-1)}] 
	\node at (0,0) {$,$};
\end{scope}
\end{tikzpicture}
\hspace{30pt}
\begin{tikzpicture}[string, yscale=0.4, xscale=0.8]
\begin{scope}[shift={(-2.5,0)}]
\begin{pgfonlayer}{bg}
\begin{scope}[shift={(0.5,0.98)}]
	\path[fill, color=gray!5] (-1.5,-1) -- (0.5,-1) to [out=90, in=-60] (0,0) -- (-0.5,1) -- (-1.5,1) -- cycle;
	\path[fill, color=cyan!100] (-0.5,1) to (0,0) to [out=-60, in=90] (0.5,-1) -- (1,-1) -- (1,1) --  cycle;
\end{scope}
\begin{scope}[shift={(-0.5,-0.98)}]
	\path[fill, color=gray!5] (-0.5,-1) rectangle (1.5,1);
	\path[fill, color=cyan!100] (1.5,-1) rectangle (2,1);
\end{scope}
\end{pgfonlayer}
\begin{pgfonlayer}{mid}
\begin{scope}[shift={(0.5,1)}]
	\draw[edge] (-0.5, -1) to (0,0);
	\draw[edge] (-0.5, 1) to (0,0);
	\draw[edge, out=90, in=-60] (0.5, -1) to (0,0);
	\node[dotdark] at (0,0) {};
\end{scope}
\begin{scope}[shift={(-0.5,-1)}]
	\draw[edge] (0, 0) to (0.5,1);
	\draw[edge] (1.5,-1) to (1.5,1);
	\node[dot] at (0,0) {};
\end{scope}
\end{pgfonlayer}
\end{scope}
\begin{scope}
	\draw[3c1] (-0.8,0) to (0.8,0);
	\draw[3c2] (-0.8,0) to (0.8,0);
	\draw[3c3] (-0.8,0) to (0.8,0);
\end{scope}
\begin{scope}[shift={(2.5,0)}]
\begin{pgfonlayer}{bg}
\begin{scope}[shift={(0.5,0.98)}]
	\path[fill, color=gray!5] (-2,-1) rectangle (-1,1);
	\path[fill, color=cyan!100] (-1,-1) rectangle (0,1);
\end{scope}
\begin{scope}[shift={(0.5,-0.98)}]
	\path[fill, color=gray!5] (-2,-1) rectangle (-1,1);
	\path[fill, color=cyan!100] (-1,-1) rectangle (0,1);
\end{scope}
\end{pgfonlayer}
\begin{pgfonlayer}{mid}
\begin{scope}[shift={(0.5,1)}]
	\draw[edge] (-1,-1) to (-1,1);
\end{scope}
\begin{scope}[shift={(0.5,-1)}]
	\draw[edge] (-1,-1) to (-1,1);
\end{scope}
\end{pgfonlayer}
\end{scope}
\end{tikzpicture}
\end{equation*}
Dually, the past cone captures a notion of \emph{right co-action}, and reversing cells gives all the usual mirror variants. Generalising how we obtained (\ref{eq:frobenius}), we can also construct theories of objects with a compatible left (co-)action of $X$ and right (co-)action of $Y$, as a quotient of $I \otimes X \otimes Y$.

\begin{remark} It may be worth observing that, in the compositional style of building algebraic theories, the ``simplest'' theories are the ones with the least artificial identifications of cells --- hence, the ones with the \emph{most} cells, or most ``colours''. Thus, for instance, the theory of semigroups arises more naturally as a quotient of the theory of actions of non-associative operations; or the theory of Frobenius algebras as a quotient of the theory of objects with a compatible left action of a monoid and right co-action of a comonoid. This is in contrast with direct, symbolic presentations, where one is led to consider theories with a small signature as simpler.
\end{remark}

\section{Smash products: commutativity, bialgebras} \label{sec:smash}
When $X$ and $Y$ are two PRO-like computads, hence presentations of algebraic theories in the standard, narrower sense, we may want to be able to compose them in order to obtain another. Starting from the tensor product, if $a_X$, $a_Y$ are the basic object cells of $X \otimes Y$, the obvious candidate for the new basic cell is $a_X \otimes a_Y$; thus, we need to quotient out the $a_X \otimes *$ and $* \otimes a_Y$. This comes naturally if we work with \emph{pointed spaces}.
\begin{dfn}
A \emph{pointed computad} $(X, *_X)$ is a computad $X$ with a distinguished $0$-cell $*_X$, its \emph{basepoint}. A map of pointed computads $f: (X, *_X) \to (Y, *_Y)$ is a map of computads with $f(*_X) = *_Y$.

Given two pointed computads $(X, *_X)$, $(Y, *_Y)$, their \emph{wedge sum} is the pointed computad $(X \lor Y, *)$, where $X \lor Y := (X \oplus Y)/(\{*_X\} \oplus \{*_Y\})$, and $*$ is the identification of $*_X$ and $*_Y$.

There is an inclusion of computads $X \lor Y \subseteq X \otimes Y$, given by $X \mapsto X \otimes \{*_Y\}$, $Y \mapsto \{*_X\} \otimes Y$. The \emph{smash product} of $X$ and $Y$ is the pointed computad $X \land Y := (X \otimes Y)/(X \lor Y)$, with the image of $X \lor Y$ through the quotient as basepoint.
\end{dfn}

\begin{exm}
The reduced cylinder of $(X, *_X)$, mentioned in the previous section, can be described as the smash product $(X, *_X) \land (I + 1, *_1)$.
\end{exm}

Much like their analogues in topology, the wedge sum and smash product define monoidal structures on the category $\mathbf{Cpt}_*$ of pointed computads, with $(1, *)$ and $(1 + 1, *)$, respectively, as units. We will sometimes leave the basepoint implicit.

Let $S^1$ be the \emph{oriented $1$-sphere}, that is, the computad
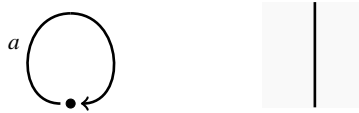
\begin{equation*}
\begin{tikzpicture}
\begin{pgfonlayer}{mid}
	\node[0c] (0) at (0,-1.2) {};
	\draw[edge, shorten <=2pt] (0) to [out=180, in=180, looseness=1.5] node[auto] {$a$} (0,0);
	\draw[edge, ->, shorten >=2pt] (0,0) to [out=0, in=0, looseness=1.5] (0);
\end{pgfonlayer} 
\end{tikzpicture}
\hspace{50pt}
\begin{tikzpicture}[string, scale=0.7]
\begin{pgfonlayer}{bg}
	\fill[color=gray!5] (-1,-1) rectangle (1,1);
\end{pgfonlayer}
\begin{pgfonlayer}{mid}
	\draw[edge] (0,-1) to (0,1);
\end{pgfonlayer}
\end{tikzpicture}
\end{equation*}
with its unique 0-cell as basepoint.
\begin{dfn}
The \emph{reduced suspension} of a pointed computad $(X, *)$ is the computad $\Sigma X := X \land S^1$.
\end{dfn}
Write $\sigma \land \tau$ for the image of $\sigma \otimes \tau$ through the quotient. Generators of $\Sigma X$ of dimension $d > 0$ are in bijection with generators of $X$, by the assignment $\sigma \mapsto \sigma \land a$, and $\partial^\alpha_n(\sigma \land a) = \partial^\alpha_n \sigma \land a$, for all numbers $n$, and $\alpha \in \{+,-\}$.

The net effect of the suspension on a PRO-like computad $X$, thus, is just to raise the dimension of each cell --- and, in particular, the dimension of the ``objects'' --- by 1. This is useful when one needs to compare theories of different basic dimensionality: the multiplication may be represented by a 2-cell in a theory of monoids $M$, and by a 3-cell in a theory of commutative monoids $M_\mathrm{comm}$, yet we want to be able to identify the two; the solution is to include $\Sigma M$, rather than $M$, in $M_\mathrm{comm}$.

Next, let us consider the smash product of the theory of monoids with itself, $M \land M$. The 3-cells $\mu \land a$ and $\eta \land a$ are suspensions of a multiplication and unit. In the diagrammatic representation, their arity as operations is a reflection of the sliding moves changing the number of intersections between diagrams from one copy of $M$ and the other; these are the only 2-cells of $M \otimes M$ that survive the quotient.
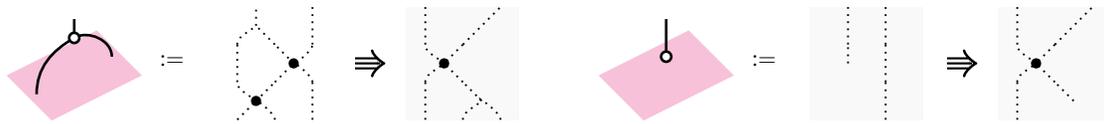
\begin{equation*}
\begin{tikzpicture} 
\begin{scope}[shift={(-1.5,0)}, scale=0.5]
\begin{pgfonlayer}{bg}
	\path[fill, color=magenta!30] (-1.8,0) -- (0.6,1.2) -- (1.8,0) -- (-0.6,-1.2) -- cycle;
\end{pgfonlayer}
\begin{pgfonlayer}{mid}
	\draw[edge] (0,1) to (0,1.5);
	\draw[edge, out=90, in=-150] (-1,-0.5) to (0,1);
	\draw[edge, out=90, in=30] (1,0.5) to (0,1);
\end{pgfonlayer}
\begin{pgfonlayer}{fg}
	\node[dot] at (0,1) {};
\end{pgfonlayer}
\end{scope}
\begin{scope}[shift={(-0.2,0.2)}] 
	\node at (0,0) {$:=$};
\end{scope}
\end{tikzpicture}
\hspace{10pt}
\begin{tikzpicture}[string, yscale=0.25, xscale=-0.5] 
\begin{scope}[shift={(-2.5,0)}]
\begin{pgfonlayer}{bg}
	\path[fill, color=gray!5] (-1.5,-3) rectangle (1.5,3);
\end{pgfonlayer}
\begin{pgfonlayer}{mid}
\begin{scope}[shift={(-0.5,2)}]
	\draw[edgedotdark] (0.5, -1) to (-0.5,1);
	\draw[edgedotdark] (1.5,-1) to (1.5,1);
\end{scope}
\begin{scope}[shift={(0.5,0)}]
	\draw[edgedotdark] (-0.5, -1) to (0,0);
	\draw[edgedotdark] (-0.5, 1) to (0,0);
	\draw[edgedotdark, out=60, in=-90] (0, 0) to (0.5,1);
	\draw[edgedotdark, out=90, in=-60] (0.5, -1) to (0,0);
	\node[dotdark] at (0,0) {};
\end{scope}
\begin{scope}[shift={(-0.5,-2)}]
	\draw[edgedotdark, out=90, in=60] (-0.5, -1) to (0,0);
	\draw[edgedotdark] (0, 0) to (0.5,1);
	\draw[edgedotdark, out=90, in=-60] (0.5, -1) to (0,0);
	\draw[edgedotdark] (1.5,-1) to (1.5,1);
\end{scope}
\end{pgfonlayer}
\end{scope}
\begin{scope}
	\draw[3c1] (1,0) to (-1,0);
	\draw[3c2] (1,0) to (-1,0);
	\draw[3c3] (1,0) to (-1,0);
\end{scope}
\begin{scope}[shift={(2.5,0)}]
\begin{pgfonlayer}{bg}
\begin{pgfonlayer}{bg}
	\path[fill, color=gray!5] (-1.5,-3) rectangle (1.5,3);
\end{pgfonlayer}
\end{pgfonlayer}
\begin{pgfonlayer}{mid}
\begin{scope}[shift={(0.5,2)}]
	\draw[edgedotdark] (-0.5, -1) to (0,0);
	\draw[edgedotdark] (0, 0) to (0,1);
	\draw[edgedotdark, out=90, in=-60] (0.5, -1) to (0,0);
	\draw[edgedotdark] (-1.5,-1) to (-1.5,1);
\end{scope}
\begin{scope}[shift={(-0.5,0)}]
	\draw[edgedotdark, out=90, in=60] (-0.5, -1) to (0,0);
	\draw[edgedotdark] (0, 0) to (0.5,1);
	\draw[edgedotdark] (0.5, -1) to (0,0);
	\draw[edgedotdark, out=-60, in=-90] (0, 0) to (-0.5,1);
	\draw[edgedotdark] (1.5,-1) to (1.5,1);
	\node[dotdark] at (0,0) {};
\end{scope}
\begin{scope}[shift={(0.5,-2)}]
	\draw[edgedotdark, out=90, in=60] (-0.5, -1) to (0,0);
	\draw[edgedotdark, out=60, in=-90] (0, 0) to (0.5,1);
	\draw[edgedotdark] (0.5, -1) to (0,0);
	\draw[edgedotdark] (-1.5,-1) to (-1.5,1);
	\draw[edgedotdark] (0, 0) to (-0.5,1);
	\node[dotdark] at (0,0) {};
\end{scope}
\end{pgfonlayer}
\end{scope}
\end{tikzpicture}
\hspace{30pt}
\begin{tikzpicture} 
\begin{scope}[shift={(-1.5,0)}, scale=0.5]
\begin{pgfonlayer}{bg}
	\path[fill, color=magenta!30] (-1.8,0) -- (0.6,1.2) -- (1.8,0) -- (-0.6,-1.2) -- cycle;
\end{pgfonlayer}
\begin{pgfonlayer}{mid}
	\draw[edge] (0,0.5) to (0,1.5);
\end{pgfonlayer}
\begin{pgfonlayer}{fg}
	\node[dot] at (0,0.5) {};
\end{pgfonlayer}
\end{scope}
\begin{scope}[shift={(-0.2,0.2)}] 
	\node at (0,0) {$:=$};
\end{scope}
\end{tikzpicture}
\hspace{10pt}
\begin{tikzpicture}[string, yscale=0.25, xscale=-0.5] 
\begin{scope}[shift={(-2.5,0)}]
\begin{pgfonlayer}{bg}
	\path[fill, color=gray!5] (-1.5,-3) rectangle (1.5,3);
\end{pgfonlayer}
\begin{pgfonlayer}{mid}
\begin{scope}[shift={(-0.5,2)}]
	\draw[edgedotdark] (0.5, -1) to (-0.5,1);
	\draw[edgedotdark] (1.5,-1) to (1.5,1);
\end{scope}
\begin{scope}[shift={(0.5,0)}]
	\draw[edgedotdark] (-0.5, -1) to (0,0);
	\draw[edgedotdark] (-0.5, 1) to (0,0);
	\draw[edgedotdark, out=60, in=-90] (0, 0) to (0.5,1);
	\draw[edgedotdark, out=90, in=-60] (0.5, -1) to (0,0);
	\node[dotdark] at (0,0) {};
\end{scope}
\begin{scope}[shift={(-0.5,-2)}]
	\draw[edgedotdark] (0, 0) to (0.5,1);
	\draw[edgedotdark] (1.5,-1) to (1.5,1);
\end{scope}
\end{pgfonlayer}
\end{scope}
\begin{scope}
	\draw[3c1] (1,0) to (-1,0);
	\draw[3c2] (1,0) to (-1,0);
	\draw[3c3] (1,0) to (-1,0);
\end{scope}
\begin{scope}[shift={(2.5,0)}]
\begin{pgfonlayer}{bg}
	\path[fill, color=gray!5] (-1.5,-3) rectangle (1.5,3);
\end{pgfonlayer}
\begin{pgfonlayer}{mid}
\begin{scope}[shift={(0.5,2)}]
	\draw[edgedotdark] (-1,-1) to (-1,1);
	\draw[edgedotdark] (0,1) to (0,-1);
\end{scope}
\begin{scope}[shift={(0.5,0)}]
	\draw[edgedotdark] (-1,-1) to (-1,1);
	\draw[edgedotdark] (0,0) to (0,1);
\end{scope}
\begin{scope}[shift={(0.5,-2)}]
	\draw[edgedotdark] (-1,-1) to (-1,1);
\end{scope}
\end{pgfonlayer}
\end{scope}
\end{tikzpicture}
\end{equation*}
On the other hand, since $a$ is an odd-dimensional cell, $a \land \mu$ and $a \land \eta$ are suspensions of a \emph{co}-multiplication and \emph{co}-unit.
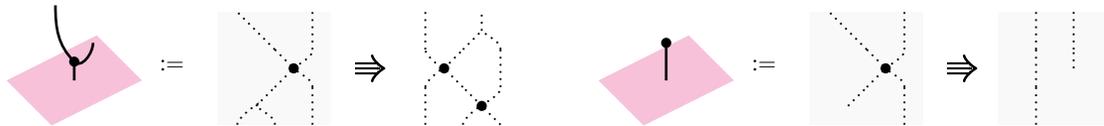
\begin{equation*}
\begin{tikzpicture} 
\begin{scope}[shift={(-1.5,0)}, scale=0.5]
\begin{pgfonlayer}{bg}
	\path[fill, color=magenta!30] (-1.8,0) -- (0.6,1.2) -- (1.8,0) -- (-0.6,-1.2) -- cycle;
\end{pgfonlayer}
\begin{pgfonlayer}{mid}
	\draw[edge] (0,0) to (0,0.5);
	\draw[edge, out=-45, in=-90] (0,0.5) to (0.5,1);
	\draw[edge, out=135, in=-90] (0,0.5) to (-0.5, 2);
\end{pgfonlayer}
\begin{pgfonlayer}{fg}
	\node[dotdark] at (0,0.5) {};
\end{pgfonlayer}
\end{scope}
\begin{scope}[shift={(-0.2,0.2)}] 
	\node at (0,0) {$:=$};
\end{scope}
\end{tikzpicture}
\hspace{10pt}
\begin{tikzpicture}[string, yscale=0.25, xscale=0.5] 
\begin{scope}[shift={(-2.5,0)}]
\begin{pgfonlayer}{bg}
	\path[fill, color=gray!5] (-1.5,-3) rectangle (1.5,3);
\end{pgfonlayer}
\begin{pgfonlayer}{mid}
\begin{scope}[shift={(-0.5,2)}]
	\draw[edgedotdark] (0.5, -1) to (-0.5,1);
	\draw[edgedotdark] (1.5,-1) to (1.5,1);
\end{scope}
\begin{scope}[shift={(0.5,0)}]
	\draw[edgedotdark] (-0.5, -1) to (0,0);
	\draw[edgedotdark] (-0.5, 1) to (0,0);
	\draw[edgedotdark, out=60, in=-90] (0, 0) to (0.5,1);
	\draw[edgedotdark, out=90, in=-60] (0.5, -1) to (0,0);
	\node[dotdark] at (0,0) {};
\end{scope}
\begin{scope}[shift={(-0.5,-2)}]
	\draw[edgedotdark, out=90, in=60] (-0.5, -1) to (0,0);
	\draw[edgedotdark] (0, 0) to (0.5,1);
	\draw[edgedotdark, out=90, in=-60] (0.5, -1) to (0,0);
	\draw[edgedotdark] (1.5,-1) to (1.5,1);
\end{scope}
\end{pgfonlayer}
\end{scope}
\begin{scope}
	\draw[3c1] (-1,0) to (1,0);
	\draw[3c2] (-1,0) to (1,0);
	\draw[3c3] (-1,0) to (1,0);
\end{scope}
\begin{scope}[shift={(2.5,0)}]
\begin{pgfonlayer}{bg}
\begin{pgfonlayer}{bg}
	\path[fill, color=gray!5] (-1.5,-3) rectangle (1.5,3);
\end{pgfonlayer}
\end{pgfonlayer}
\begin{pgfonlayer}{mid}
\begin{scope}[shift={(0.5,2)}]
	\draw[edgedotdark] (-0.5, -1) to (0,0);
	\draw[edgedotdark] (0, 0) to (0,1);
	\draw[edgedotdark, out=90, in=-60] (0.5, -1) to (0,0);
	\draw[edgedotdark] (-1.5,-1) to (-1.5,1);
\end{scope}
\begin{scope}[shift={(-0.5,0)}]
	\draw[edgedotdark, out=90, in=60] (-0.5, -1) to (0,0);
	\draw[edgedotdark] (0, 0) to (0.5,1);
	\draw[edgedotdark] (0.5, -1) to (0,0);
	\draw[edgedotdark, out=-60, in=-90] (0, 0) to (-0.5,1);
	\draw[edgedotdark] (1.5,-1) to (1.5,1);
	\node[dotdark] at (0,0) {};
\end{scope}
\begin{scope}[shift={(0.5,-2)}]
	\draw[edgedotdark, out=90, in=60] (-0.5, -1) to (0,0);
	\draw[edgedotdark, out=60, in=-90] (0, 0) to (0.5,1);
	\draw[edgedotdark] (0.5, -1) to (0,0);
	\draw[edgedotdark] (-1.5,-1) to (-1.5,1);
	\draw[edgedotdark] (0, 0) to (-0.5,1);
	\node[dotdark] at (0,0) {};
\end{scope}
\end{pgfonlayer}
\end{scope}
\end{tikzpicture}
\hspace{30pt}
\begin{tikzpicture} 
\begin{scope}[shift={(-1.5,0)}, scale=0.5]
\begin{pgfonlayer}{bg}
	\path[fill, color=magenta!30] (-1.8,0) -- (0.6,1.2) -- (1.8,0) -- (-0.6,-1.2) -- cycle;
\end{pgfonlayer}
\begin{pgfonlayer}{mid}
	\draw[edge] (0,0) to (0,1);
\end{pgfonlayer}
\begin{pgfonlayer}{fg}
	\node[dotdark] at (0,1) {};
\end{pgfonlayer}
\end{scope}
\begin{scope}[shift={(-0.2,0.2)}] 
	\node at (0,0) {$:=$};
\end{scope}
\end{tikzpicture}
\hspace{10pt}
\begin{tikzpicture}[string, yscale=0.25, xscale=0.5] 
\begin{scope}[shift={(-2.5,0)}]
\begin{pgfonlayer}{bg}
	\path[fill, color=gray!5] (-1.5,-3) rectangle (1.5,3);
\end{pgfonlayer}
\begin{pgfonlayer}{mid}
\begin{scope}[shift={(-0.5,2)}]
	\draw[edgedotdark] (0.5, -1) to (-0.5,1);
	\draw[edgedotdark] (1.5,-1) to (1.5,1);
\end{scope}
\begin{scope}[shift={(0.5,0)}]
	\draw[edgedotdark] (-0.5, -1) to (0,0);
	\draw[edgedotdark] (-0.5, 1) to (0,0);
	\draw[edgedotdark, out=60, in=-90] (0, 0) to (0.5,1);
	\draw[edgedotdark, out=90, in=-60] (0.5, -1) to (0,0);
	\node[dotdark] at (0,0) {};
\end{scope}
\begin{scope}[shift={(-0.5,-2)}]
	\draw[edgedotdark] (0, 0) to (0.5,1);
	\draw[edgedotdark] (1.5,-1) to (1.5,1);
\end{scope}
\end{pgfonlayer}
\end{scope}
\begin{scope}
	\draw[3c1] (-1,0) to (1,0);
	\draw[3c2] (-1,0) to (1,0);
	\draw[3c3] (-1,0) to (1,0);
\end{scope}
\begin{scope}[shift={(2.5,0)}]
\begin{pgfonlayer}{bg}
	\path[fill, color=gray!5] (-1.5,-3) rectangle (1.5,3);
\end{pgfonlayer}
\begin{pgfonlayer}{mid}
\begin{scope}[shift={(0.5,2)}]
	\draw[edgedotdark] (-1,-1) to (-1,1);
	\draw[edgedotdark] (0,1) to (0,-1);
\end{scope}
\begin{scope}[shift={(0.5,0)}]
	\draw[edgedotdark] (-1,-1) to (-1,1);
	\draw[edgedotdark] (0,0) to (0,1);
\end{scope}
\begin{scope}[shift={(0.5,-2)}]
	\draw[edgedotdark] (-1,-1) to (-1,1);
\end{scope}
\end{pgfonlayer}
\end{scope}
\end{tikzpicture}
\end{equation*}
In the tensor product $M \otimes M$, the 4-cell $\mu \otimes \mu$ mediates between two ways of sliding diagrams past each other. In the smash product $M \land M$, this becomes the \emph{bialgebra law} between multiplication and comultiplication. 
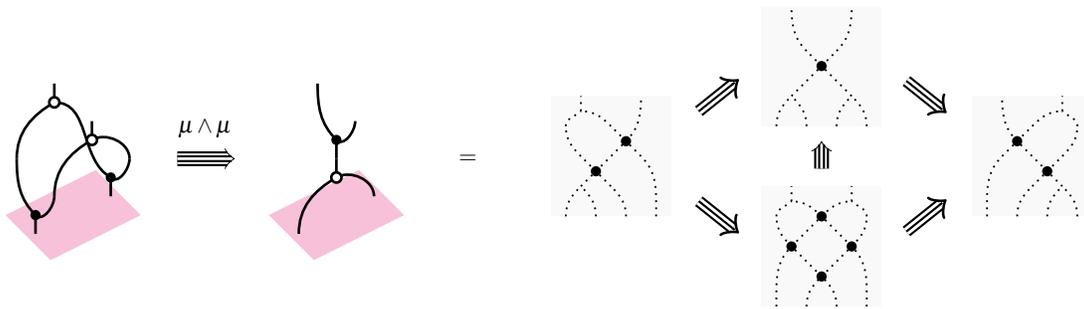
\begin{equation*}
\begin{tikzpicture}[scale=0.5, baseline=-35pt] 
\begin{scope}[shift={(-3.5,0)}] 
\begin{pgfonlayer}{bg}
	\path[fill, color=magenta!30] (-1.8,0) -- (0.6,1.2) -- (1.8,0) -- (-0.6,-1.2) -- cycle;
\end{pgfonlayer}
\begin{scope}[shift={(-1,-0.5)}] 
\begin{pgfonlayer}{mid}
	\draw[edge] (0,0) to (0,0.5);
	\draw[edge, out=-45, in=-90] (0,0.5) to (0.5,1);
	\draw[edge, out=135, in=-90] (0,0.5) to (-0.5, 2);
\end{pgfonlayer}
\begin{pgfonlayer}{fg}
	\node[dotdark] at (0,0.5) {};
\end{pgfonlayer}
\end{scope}
\begin{scope}[shift={(1,0.5)}] 
\begin{pgfonlayer}{mid}
	\draw[edge] (0,0) to (0,0.5);
	\draw[edge, out=-45, in=-90] (0,0.5) to (0.5,1);
	\draw[edge, out=150, in=0] (0,0.5) to (-1.5, 2.5);
\end{pgfonlayer}
\begin{pgfonlayer}{fg}
	\node[dotdark] at (0,0.5) {};
\end{pgfonlayer}
\end{scope}
\begin{scope}[shift={(0.5,1)}] 
\begin{pgfonlayer}{mid}
	\draw[edge] (0,1) to (0,1.5);
	\draw[edge, out=90, in=-150] (-1,-0.5) to (0,1);
	\draw[edge, out=90, in=30] (1,0.5) to (0,1);
\end{pgfonlayer}
\begin{pgfonlayer}{fg}
	\node[dot] at (0,1) {};
\end{pgfonlayer}
\end{scope}
\begin{scope}[shift={(-0.5,2)}] 
\begin{pgfonlayer}{mid}
	\draw[edge] (0,1) to (0,1.5);
	\draw[edge, out=90, in=-150] (-1,-0.5) to (0,1);
\end{pgfonlayer}
\begin{pgfonlayer}{fg}
	\node[dot] at (0,1) {};
\end{pgfonlayer}
\end{scope}
\end{scope}
\begin{scope}[shift={(3.5,0)}] 
\begin{pgfonlayer}{bg}
	\path[fill, color=magenta!30] (-1.8,0) -- (0.6,1.2) -- (1.8,0) -- (-0.6,-1.2) -- cycle;
\end{pgfonlayer}
\begin{scope} 
\begin{pgfonlayer}{mid}
	\draw[edge] (0,1) to (0,1.5);
	\draw[edge, out=90, in=-150] (-1,-0.5) to (0,1);
	\draw[edge, out=90, in=30] (1,0.5) to (0,1);
\end{pgfonlayer}
\begin{pgfonlayer}{fg}
	\node[dot] at (0,1) {};
\end{pgfonlayer}
\end{scope}
\begin{scope}[shift={(0,1.5)}] 
\begin{pgfonlayer}{mid}
	\draw[edge] (0,0) to (0,0.5);
	\draw[edge, out=-45, in=-90] (0,0.5) to (0.5,1);
	\draw[edge, out=135, in=-90] (0,0.5) to (-0.5, 2);
\end{pgfonlayer}
\begin{pgfonlayer}{fg}
	\node[dotdark] at (0,0.5) {};
\end{pgfonlayer}
\end{scope}
\end{scope}
\begin{scope}[shift={(0,1.5)}] 
	\draw[4c1] (-0.8,0) to (0.8,0);
	\draw[4c2] (-0.8,0) to (0.8,0);
	\node at (0,0.8) {$\mu \land \mu$};
\end{scope}
\begin{scope}[shift={(7,1.5)}] 
	\node at (0,0) {$=$};
\end{scope}
\end{tikzpicture}
\hspace{25pt}
\begin{tikzpicture}[scale=0.4]
\begin{scope}[shift={(-7,0)}] 
\begin{pgfonlayer}{bg}
	\path[fill, color=gray!5] (-2,-2) rectangle (2,2);
\end{pgfonlayer}
\begin{pgfonlayer}{mid}
	\draw[edgedotdark, out=90, in=-135] (-1.5,-2) to (-1, -1); 
	\draw[edgedotdark, out=90, in=-45] (-0.5,-2) to (-1, -1);
	\draw[edgedotdark] (-1,-1) to (0.5, 0.5);
	\draw[edgedotdark, in=-90] (0.5,0.5) to (1, 2);
	\draw[edgedotdark, out=90, in=-45] (0.5,-2) to (0,-1); 
	\draw[edgedotdark] (0,-1) to (-1,0);
	\draw[edgedotdark, out=90, in=-45] (1.5,-2) to (1,0);
	\draw[edgedotdark] (1,0) to (0,1);
	\draw[edgedotdark, out=135, in=-90] (-1,0) to (-1.5,1);
	\draw[edgedotdark, out=90, in=-135] (-1.5,1) to (-1,1.5);
	\draw[edgedotdark, out=135, in=0] (0,1) to (-1,1.5);
	\draw[edgedotdark] (-1,1.5) to (-1,2);
\end{pgfonlayer}
\begin{pgfonlayer}{fg}
	\node[dotdark] at (-0.5,-0.5) {};
	\node[dotdark] at (0.5,0.5) {};
\end{pgfonlayer}
\end{scope}
\begin{scope}[shift={(0,-3)}] 
\begin{pgfonlayer}{bg}
	\path[fill, color=gray!5] (-2,-2) rectangle (2,2);
\end{pgfonlayer}
\begin{pgfonlayer}{mid}
	\draw[edgedotdark, out=90, in=-135] (-0.5,-2) to (0,-1); 
	\draw[edgedotdark] (0,-1) to (1,0);
	\draw[edgedotdark, out=90, in=-135] (-1.5,-2) to (-1,0);
	\draw[edgedotdark] (-1,0) to (0,1);
	\draw[edgedotdark, out=45, in=-90] (1,0) to (1.5,1);
	\draw[edgedotdark, out=90, in=-45] (1.5,1) to (1,1.5);
	\draw[edgedotdark, out=45, in=180] (0,1) to (1,1.5);
	\draw[edgedotdark] (1,1.5) to (1,2);
	\draw[edgedotdark, out=90, in=-45] (0.5,-2) to (0,-1); 
	\draw[edgedotdark] (0,-1) to (-1,0);
	\draw[edgedotdark, out=90, in=-45] (1.5,-2) to (1,0);
	\draw[edgedotdark] (1,0) to (0,1);
	\draw[edgedotdark, out=135, in=-90] (-1,0) to (-1.5,1);
	\draw[edgedotdark, out=90, in=-135] (-1.5,1) to (-1,1.5);
	\draw[edgedotdark, out=135, in=0] (0,1) to (-1,1.5);
	\draw[edgedotdark] (-1,1.5) to (-1,2);
\end{pgfonlayer}
\begin{pgfonlayer}{fg}
	\node[dotdark] at (-1,0) {};
	\node[dotdark] at (0,-1) {};
	\node[dotdark] at (1,0) {};
	\node[dotdark] at (0,1) {};
\end{pgfonlayer}
\end{scope}
\begin{scope}[shift={(7,0)}, xscale=-1] 
\begin{pgfonlayer}{bg}
	\path[fill, color=gray!5] (-2,-2) rectangle (2,2);
\end{pgfonlayer}
\begin{pgfonlayer}{mid}
	\draw[edgedotdark, out=90, in=-135] (-1.5,-2) to (-1, -1); 
	\draw[edgedotdark, out=90, in=-45] (-0.5,-2) to (-1, -1);
	\draw[edgedotdark] (-1,-1) to (0.5, 0.5);
	\draw[edgedotdark, in=-90] (0.5,0.5) to (1, 2);
	\draw[edgedotdark, out=90, in=-45] (0.5,-2) to (0,-1); 
	\draw[edgedotdark] (0,-1) to (-1,0);
	\draw[edgedotdark, out=90, in=-45] (1.5,-2) to (1,0);
	\draw[edgedotdark] (1,0) to (0,1);
	\draw[edgedotdark, out=135, in=-90] (-1,0) to (-1.5,1);
	\draw[edgedotdark, out=90, in=-135] (-1.5,1) to (-1,1.5);
	\draw[edgedotdark, out=135, in=0] (0,1) to (-1,1.5);
	\draw[edgedotdark] (-1,1.5) to (-1,2);
\end{pgfonlayer}
\begin{pgfonlayer}{fg}
	\node[dotdark] at (-0.5,-0.5) {};
	\node[dotdark] at (0.5,0.5) {};
\end{pgfonlayer}
\end{scope}
\begin{scope}[shift={(0,3)}] 
\begin{pgfonlayer}{bg}
	\path[fill, color=gray!5] (-2,-2) rectangle (2,2);
\end{pgfonlayer}
\begin{pgfonlayer}{mid}
	\draw[edgedotdark, out=90, in=-135] (-1.5,-2) to (-1, -1); 
	\draw[edgedotdark, out=90, in=-45] (-0.5,-2) to (-1, -1);
	\draw[edgedotdark] (-1,-1) to (0.5, 0.5);
	\draw[edgedotdark, in=-90] (0.5,0.5) to (1, 2);
	\draw[edgedotdark, out=90, in=-45] (1.5,-2) to (1, -1); 
	\draw[edgedotdark, out=90, in=-135] (0.5,-2) to (1, -1);
	\draw[edgedotdark] (1,-1) to (-0.5, 0.5);
	\draw[edgedotdark, out=135, in=-90] (-0.5,0.5) to (-1, 2);
\end{pgfonlayer}
\begin{pgfonlayer}{fg}
	\node[dotdark] at (0,0) {};
\end{pgfonlayer}
\end{scope}
\begin{scope}[shift={(-3.2,3)}] 
	\draw[3c1] (-1.5,-2) to (1,0);
	\draw[3c2] (-1.5,-2) to (1,0);
	\draw[3c3] (-1.5,-2) to (1,0);
\end{scope}
\begin{scope}[shift={(3.2,3)}] 
	\draw[3c1] (-1,0) to (1.5,-2);
	\draw[3c2] (-1,0) to (1.5,-2);
	\draw[3c3] (-1,0) to (1.5,-2);
\end{scope}
\begin{scope}[shift={(-3.2,-3)}] 
	\draw[3c1] (-1.5,2) to (1,0);
	\draw[3c2] (-1.5,2) to (1,0);
	\draw[3c3] (-1.5,2) to (1,0);
\end{scope}
\begin{scope}[shift={(3.2,-3)}] 
	\draw[3c1] (-1,0) to (1.5,2);
	\draw[3c2] (-1,0) to (1.5,2);
	\draw[3c3] (-1,0) to (1.5,2);
\end{scope}
\begin{scope}[shift={(0,0)}] 
	\draw[4c1] (0,-0.5) to (0,0.5);
	\draw[4c2] (0,-0.5) to (0,0.5);
\end{scope}
\end{tikzpicture}
\end{equation*}
The two 3-cells in the input 3-border should really be decomposed in two more, sliding each multiplication first past one strand, then past the other. The 4-cells $\mu \land \eta$, $\eta \land \mu$ and $\eta \land \eta$ give the remaining bialgebra equations.
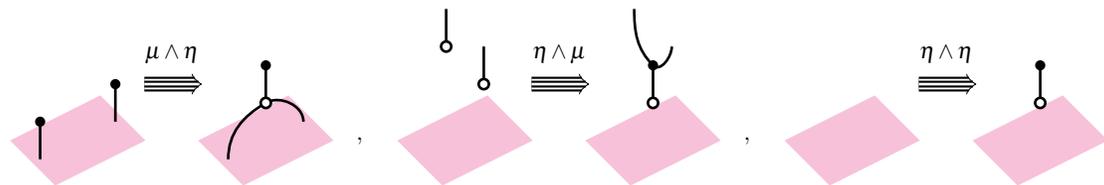
\begin{equation*}
\begin{tikzpicture}[scale=0.5] 
\begin{scope}[shift={(-2.5,0)}] 
\begin{pgfonlayer}{bg}
	\path[fill, color=magenta!30] (-1.8,0) -- (0.6,1.2) -- (1.8,0) -- (-0.6,-1.2) -- cycle;
\end{pgfonlayer}
\begin{scope}[shift={(-1,-0.5)}] 
\begin{pgfonlayer}{mid}
	\draw[edge] (0,0) to (0,1);
\end{pgfonlayer}
\begin{pgfonlayer}{fg}
	\node[dotdark] at (0,1) {};
\end{pgfonlayer}
\end{scope}
\begin{scope}[shift={(1,0.5)}] 
\begin{pgfonlayer}{mid}
	\draw[edge] (0,0) to (0,1);
\end{pgfonlayer}
\begin{pgfonlayer}{fg}
	\node[dotdark] at (0,1) {};
\end{pgfonlayer}
\end{scope}
\end{scope}
\begin{scope}[shift={(2.5,0)}] 
\begin{pgfonlayer}{bg}
	\path[fill, color=magenta!30] (-1.8,0) -- (0.6,1.2) -- (1.8,0) -- (-0.6,-1.2) -- cycle;
\end{pgfonlayer}
\begin{scope} 
\begin{pgfonlayer}{mid}
	\draw[edge] (0,1) to (0,1.5);
	\draw[edge, out=90, in=-150] (-1,-0.5) to (0,1);
	\draw[edge, out=90, in=30] (1,0.5) to (0,1);
\end{pgfonlayer}
\begin{pgfonlayer}{fg}
	\node[dot] at (0,1) {};
\end{pgfonlayer}
\end{scope}
\begin{scope}[shift={(0,1.5)}] 
\begin{pgfonlayer}{mid}
	\draw[edge] (0,0) to (0,0.5);
\end{pgfonlayer}
\begin{pgfonlayer}{fg}
	\node[dotdark] at (0,0.5) {};
\end{pgfonlayer}
\end{scope}
\end{scope}
\begin{scope}[shift={(0,1.5)}] 
	\draw[4c1] (-0.8,0) to (0.8,0);
	\draw[4c2] (-0.8,0) to (0.8,0);
	\node at (0,0.8) {$\mu \land \eta$};
\end{scope}
\begin{scope}[shift={(5,0)}] 
	\node at (0,0) {$,$};
\end{scope}
\end{tikzpicture}
\hspace{10pt}
\begin{tikzpicture}[scale=0.5]
\begin{scope}[shift={(-2.5,0)}] 
\begin{pgfonlayer}{bg}
	\path[fill, color=magenta!30] (-1.8,0) -- (0.6,1.2) -- (1.8,0) -- (-0.6,-1.2) -- cycle;
\end{pgfonlayer}
\begin{scope}[shift={(0.5,1)}] 
\begin{pgfonlayer}{mid}
	\draw[edge] (0,0.5) to (0,1.5);
\end{pgfonlayer}
\begin{pgfonlayer}{fg}
	\node[dot] at (0,0.5) {};
\end{pgfonlayer}
\end{scope}
\begin{scope}[shift={(-0.5,2)}] 
\begin{pgfonlayer}{mid}
	\draw[edge] (0,0.5) to (0,1.5);
\end{pgfonlayer}
\begin{pgfonlayer}{fg}
	\node[dot] at (0,0.5) {};
\end{pgfonlayer}
\end{scope}
\end{scope}
\begin{scope}[shift={(2.5,0)}] 
\begin{pgfonlayer}{bg}
	\path[fill, color=magenta!30] (-1.8,0) -- (0.6,1.2) -- (1.8,0) -- (-0.6,-1.2) -- cycle;
\end{pgfonlayer}
\begin{scope} 
\begin{pgfonlayer}{mid}
	\draw[edge] (0,1) to (0,1.5);
\end{pgfonlayer}
\begin{pgfonlayer}{fg}
	\node[dot] at (0,1) {};
\end{pgfonlayer}
\end{scope}
\begin{scope}[shift={(0,1.5)}] 
\begin{pgfonlayer}{mid}
	\draw[edge] (0,0) to (0,0.5);
	\draw[edge, out=-45, in=-90] (0,0.5) to (0.5,1);
	\draw[edge, out=135, in=-90] (0,0.5) to (-0.5, 2);
\end{pgfonlayer}
\begin{pgfonlayer}{fg}
	\node[dotdark] at (0,0.5) {};
\end{pgfonlayer}
\end{scope}
\end{scope}
\begin{scope}[shift={(0,1.5)}] 
	\draw[4c1] (-0.8,0) to (0.8,0);
	\draw[4c2] (-0.8,0) to (0.8,0);
	\node at (0,0.8) {$\eta \land \mu$};
\end{scope}
\begin{scope}[shift={(5,0)}] 
	\node at (0,0) {$,$};
\end{scope}
\end{tikzpicture}
\hspace{10pt}
\begin{tikzpicture}[scale=0.5]
\begin{scope}[shift={(-2.5,0)}] 
\begin{pgfonlayer}{bg}
	\path[fill, color=magenta!30] (-1.8,0) -- (0.6,1.2) -- (1.8,0) -- (-0.6,-1.2) -- cycle;
\end{pgfonlayer}
\end{scope}
\begin{scope}[shift={(2.5,0)}] 
\begin{pgfonlayer}{bg}
	\path[fill, color=magenta!30] (-1.8,0) -- (0.6,1.2) -- (1.8,0) -- (-0.6,-1.2) -- cycle;
\end{pgfonlayer}
\begin{scope} 
\begin{pgfonlayer}{mid}
	\draw[edge] (0,1) to (0,1.5);
\end{pgfonlayer}
\begin{pgfonlayer}{fg}
	\node[dot] at (0,1) {};
\end{pgfonlayer}
\end{scope}
\begin{scope}[shift={(0,1.5)}] 
\begin{pgfonlayer}{mid}
	\draw[edge] (0,0) to (0,0.5);
\end{pgfonlayer}
\begin{pgfonlayer}{fg}
	\node[dotdark] at (0,0.5) {};
\end{pgfonlayer}
\end{scope}
\end{scope}
\begin{scope}[shift={(0,1.5)}] 
	\draw[4c1] (-0.8,0) to (0.8,0);
	\draw[4c2] (-0.8,0) to (0.8,0);
	\node at (0,0.8) {$\eta \land \eta$};
\end{scope}
\end{tikzpicture}
\end{equation*}
Therefore, $M \land M$ is a presentation of the theory of bialgebras. Of course, $M \land M$ contains other higher-dimensional cells; if $\alpha$ is the associator 3-cell, $\alpha \land a$ will be an associator for $\mu \land a$, and $a \land \alpha$ a co-associator for $a \land \mu$; while $\alpha \land \mu$ will be a higher coherence between ``associate, then use bialgebra law'' and ``use bialgebra law, then associate'', and so on. 

Let us now take a look at $M \land M^\mathrm{op}$, the smash product of the theory of monoids and of the theory of comonoids. As before, $\mu \land a^\mathrm{op}$ and $\eta \land a^\mathrm{op}$ are a (suspended) multiplication and unit cell; but now $a \land \mu^\mathrm{op}$ and $a \land \eta^\mathrm{op}$ are \emph{also} a multiplication and unit. Below, we show what $\mu \land \mu^\mathrm{op}$ and $\eta \land \eta^\mathrm{op}$ look like; the black multiplication and unit represent $a \land \mu^\mathrm{op}$ and $a \land \eta^\mathrm{op}$, respectively.
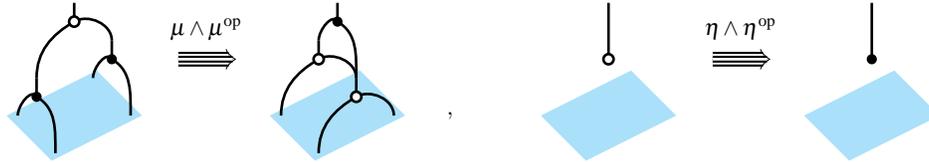
\begin{equation*}
\begin{tikzpicture}[scale=0.5]
\begin{scope}[shift={(-3.5,0)}] 
\begin{pgfonlayer}{bg}
	\path[fill, color=cyan!30] (-1.8,0) -- (0.6,1.2) -- (1.8,0) -- (-0.6,-1.2) -- cycle;
\end{pgfonlayer}
\begin{scope}[shift={(-1,-0.5)}] 
\begin{pgfonlayer}{mid}
	\draw[edge] (0,1) to (0,1.5);
	\draw[edge, out=90, in=-30] (0.5,-0.5) to (0,1);
	\draw[edge, out=90, in=120] (-0.5,0.5) to (0,1);
\end{pgfonlayer}
\begin{pgfonlayer}{fg}
	\node[dotdark] at (0,1) {};
\end{pgfonlayer}
\end{scope}
\begin{scope}[shift={(1,0.5)}] 
\begin{pgfonlayer}{mid}
	\draw[edge] (0,1) to (0,1.5);
	\draw[edge, out=90, in=-30] (0.5,-0.5) to (0,1);
	\draw[edge, out=90, in=120] (-0.5,0.5) to (0,1);
\end{pgfonlayer}
\begin{pgfonlayer}{fg}
	\node[dotdark] at (0,1) {};
\end{pgfonlayer}
\end{scope}
\begin{scope}[shift={(0,1.5)}] 
\begin{pgfonlayer}{mid}
	\draw[edge] (0,1) to (0,1.5);
	\draw[edge, out=90, in=-150] (-1,-0.5) to (0,1);
	\draw[edge, out=90, in=30] (1,0.5) to (0,1);
\end{pgfonlayer}
\begin{pgfonlayer}{fg}
	\node[dot] at (0,1) {};
\end{pgfonlayer}
\end{scope}
\end{scope}
\begin{scope}[shift={(3.5,0)}] 
\begin{pgfonlayer}{bg}
	\path[fill, color=cyan!30] (-1.8,0) -- (0.6,1.2) -- (1.8,0) -- (-0.6,-1.2) -- cycle;
\end{pgfonlayer}
\begin{scope}[shift={(0.5,-0.5)}] 
\begin{pgfonlayer}{mid}
	\draw[edge] (0,1) to (0,1.5);
	\draw[edge, out=90, in=-150] (-1,-0.5) to (0,1);
	\draw[edge, out=90, in=30] (1,0.5) to (0,1);
\end{pgfonlayer}
\begin{pgfonlayer}{fg}
	\node[dot] at (0,1) {};
\end{pgfonlayer}
\end{scope}
\begin{scope}[shift={(-0.5,0.5)}] 
\begin{pgfonlayer}{mid}
	\draw[edge] (0,1) to (0,1.5);
	\draw[edge, out=90, in=-150] (-1,-0.5) to (0,1);
	\draw[edge, out=90, in=30] (1,0.5) to (0,1);
\end{pgfonlayer}
\begin{pgfonlayer}{fg}
	\node[dot] at (0,1) {};
\end{pgfonlayer}
\end{scope}
\begin{scope}[shift={(0,1.5)}] 
\begin{pgfonlayer}{mid}
	\draw[edge] (0,1) to (0,1.5);
	\draw[edge, out=90, in=-30] (0.5,-0.5) to (0,1);
	\draw[edge, out=90, in=120] (-0.5,0.5) to (0,1);
\end{pgfonlayer}
\begin{pgfonlayer}{fg}
	\node[dotdark] at (0,1) {};
\end{pgfonlayer}
\end{scope}
\end{scope}
\begin{scope}[shift={(0,1.5)}] 
	\draw[4c1] (-0.8,0) to (0.8,0);
	\draw[4c2] (-0.8,0) to (0.8,0);
	\node at (0,0.8) {$\mu \land \mu^\mathrm{op}$};
\end{scope}
\begin{scope}[shift={(6.5,0)}] 
	\node at (0,0) {$,$};
\end{scope}
\end{tikzpicture}
\hspace{30pt}
\begin{tikzpicture}[scale=0.5]
\begin{scope}[shift={(-3.5,0)}] 
\begin{pgfonlayer}{bg}
	\path[fill, color=cyan!30] (-1.8,0) -- (0.6,1.2) -- (1.8,0) -- (-0.6,-1.2) -- cycle;
\end{pgfonlayer}
\begin{scope}[shift={(0,1.5)}] 
\begin{pgfonlayer}{mid}
	\draw[edge] (0,0) to (0,1.5);
\end{pgfonlayer}
\begin{pgfonlayer}{fg}
	\node[dot] at (0,0) {};
\end{pgfonlayer}
\end{scope}
\end{scope}
\begin{scope}[shift={(3.5,0)}] 
\begin{pgfonlayer}{bg}
	\path[fill, color=cyan!30] (-1.8,0) -- (0.6,1.2) -- (1.8,0) -- (-0.6,-1.2) -- cycle;
\end{pgfonlayer}
\begin{scope}[shift={(0,1.5)}] 
\begin{pgfonlayer}{mid}
	\draw[edge] (0,0) to (0,1.5);
\end{pgfonlayer}
\begin{pgfonlayer}{fg}
	\node[dotdark] at (0,0) {};
\end{pgfonlayer}
\end{scope}
\end{scope}
\begin{scope}[shift={(0,1.5)}] 
	\draw[4c1] (-0.8,0) to (0.8,0);
	\draw[4c2] (-0.8,0) to (0.8,0);
	\node at (0,0.8) {$\eta \land \eta^\mathrm{op}$};
\end{scope}
\end{tikzpicture}
\end{equation*}
The first one is a directed version of an \emph{interchange law}; the second one identifies the two units. By the Eckmann-Hilton argument \cite{eckmann1962group}, a pair of monoids satisfying an interchange law is equivalent to a single \emph{commutative monoid}; hence, $M \land M^\mathrm{op}$ can be seen as a presentation of the theory of commutative monoids. Dually, $M^\mathrm{op} \land M$ presents the theory of cocommutative comonoids. The defining equations of all these theories obtain an original topological interpretation in terms of intersecting, sliding diagrams.

\begin{remark}
In \cite{tubella2016subatomic}, it is suggested that directed interchange laws, of the kind just presented, may have a fundamental role in proof systems for propositional logics: what we have said about potential insights on algebraic theories coming from our approach may apply to proof systems, and their normalisation properties, as well.
\end{remark}

\begin{remark}
There is a symmetric monoidal structure on the category of symmetric operads --- the \emph{Boardman-Vogt tensor product} \cite{boardman1973homotopy} --- such that the tensor product of the operad of monoids with itself is equivalent to the operad of commutative monoids \cite{weiss2011operads}. We have not explored the connection, but by the results just presented, it seems likely that the two constructions are related, with the non-symmetric smash product of computads being the more general one.
\end{remark}

We conclude by putting all the information together in order to describe how a part of the theory of interacting bialgebras \cite{bonchi2014interacting}, the ``basic fragment'' of the ZX calculus, can be assembled in the language we developed.
\begin{enumerate}
	\item Starting from the theory of constants $K$ (either as a given, or as a 2-cube with three faces quotiented out), we can obtain theories $M$ and $M^\mathrm{op}$ of monoids and comonoids by successive cones, and identifying cells of congruent shape, as shown in Section \ref{sec:basic}
	\item By the discussion in Section \ref{sec:homo}, we know how to obtain the theory of a compatible left action of a monoid and right co-action of a comonoid by cones. A theory $F$ of Frobenius algebras results from the identification of cells of congruent shape. Also letting the comultiplication act on the left, or the multiplication co-act on the right, allows one to eliminate loops, as in the diagram below --- leading to a theory $F_s$ of special Frobenius algebras.
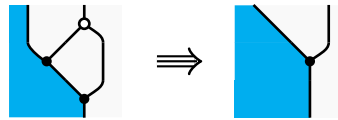
\begin{equation*}
\begin{tikzpicture}[string, yscale=0.25, xscale=0.5]
\begin{scope}[shift={(3,0)}]
\begin{pgfonlayer}{bg}
\begin{scope}[shift={(-0.5,1.98)}]
	\path[fill, color=cyan!100] (-1,-1) -- (0.5,-1) -- (-0.5, 1) -- (-1,1) -- cycle;
	\path[fill, color=gray!5] (-0.5,1) -- (0.5, -1) -- (2,-1) -- (2,1) -- cycle;
\end{scope}
\begin{scope}[shift={(0.5,0)}]
	\path[fill, color=cyan!100] (-2,-1) -- (0,-1) to (0,0) -- (-0.5,1) -- (-2,1) -- cycle;
	\path[fill, color=gray!5] (-0.5,1) -- (0,0) -- (0,-1) -- (1,-1) -- (1,1) --  cycle;
\end{scope}
\begin{scope}[shift={(-0.5,-1.98)}]
	\path[fill, color=cyan!100] (-1,-1) rectangle (1,1);
	\path[fill, color=gray!5] (1,-1) rectangle (2,1);
\end{scope}
\end{pgfonlayer}
\begin{pgfonlayer}{mid}
\begin{scope}[shift={(-0.5,2)}]
	\draw[edge] (0.5, -1) to (-0.5,1);
	\draw[edge] (1.5,-1) to (1.5,1);
\end{scope}
\begin{scope}[shift={(0.5,0)}]
	\draw[edge] (-0.5, 1) to (0,0);
	\draw[edge, out=60, in=-90] (0, 0) to (0.5,1);
	\draw[edge] (0, -1) to (0,0);
	\node[dotdark] at (0,0) {};
\end{scope}
\begin{scope}[shift={(-0.5,-2)}]
	\draw[edge] (1,-1) to (1,1);
\end{scope}
\end{pgfonlayer}
\end{scope}
\begin{scope}
	\draw[3c1] (-1.2,0) to (1.2,0);
	\draw[3c2] (-1.2,0) to (1.2,0);
	\draw[3c3] (-1.2,0) to (1.2,0);
\end{scope}
\begin{scope}[shift={(-3,0)}]
\begin{pgfonlayer}{bg}
\begin{scope}[shift={(0.5,1.98)}]
	\path[fill, color=cyan!100] (-2,-1) rectangle (-1.5,1);
	\path[fill, color=gray!5] (-1.5,-1) rectangle (1,1);
\end{scope}
\begin{scope}[shift={(-0.5,0)}]
	\path[fill, color=cyan!100] (-1,-1) -- (0.5,-1) to (0,0) to [out=-60, in= -90] (-0.5,1) -- (-1,1) -- cycle;
	\path[fill, color=gray!5] (-0.5,1) to [out=-90, in=-60] (0,0) to (0.5,-1) -- (2,-1) -- (2,1) -- cycle;
\end{scope}
\begin{scope}[shift={(0.5,-1.98)}]
	\path[fill, color=cyan!100] (-2,-1) -- (0,-1) -- (0,0) -- (-0.5,1) -- (-2,1) -- cycle;
	\path[fill, color=gray!5] (-0.5,1) -- (1,1) -- (1,-1) -- (0,-1) -- (0,0) -- cycle;
\end{scope}
\end{pgfonlayer}
\begin{pgfonlayer}{mid}
\begin{scope}[shift={(0.5,2)}]
	\draw[edge] (-0.5, -1) to (0,0);
	\draw[edge] (0, 0) to (0,1);
	\draw[edge, out=90, in=-60] (0.5, -1) to (0,0);
	\draw[edge] (-1.5,-1) to (-1.5,1);
	\node[dot] at (0,0) {};
\end{scope}
\begin{scope}[shift={(-0.5,0)}]
	\draw[edge] (0, 0) to (0.5,1);
	\draw[edge] (0.5, -1) to (0,0);
	\draw[edge, out=-60, in=-90] (0, 0) to (-0.5,1);
	\draw[edge] (1.5,-1) to (1.5,1);
	\node[dotdark] at (0,0) {};
\end{scope}
\begin{scope}[shift={(0.5,-2)}]
	\draw[edge, out=60, in=-90] (0, 0) to (0.5,1);
	\draw[edge] (0, -1) to (0,0);
	\draw[edge] (0, 0) to (-0.5,1);
	\node[dotdark] at (0,0) {};
\end{scope}
\end{pgfonlayer}
\end{scope}
\end{tikzpicture}
\end{equation*}
	\item By the results of this section, a theory of commutative, co-commutative bialgebras can be presented as a smash product of copies of $M$ and $M^\mathrm{op}$, for instance $B := M \land M \land M \land M$. 
	\item To conclude, it suffices to take two copies of $B$ and two copies of $F_s$, and identify pairs of monoids and comonoids; since objects are 4-cells in $B$ and 1-cells in $F_s$, we first need to take the iterated suspension $\Sigma^3 F_s$. Letting $A := \Sigma^3(M \lor M^\mathrm{op} \lor M \lor M^\mathrm{op})$, a presentation of the theory $IB$ of interacting bialgebras is obtained as the pushout
	\begin{equation*}
\begin{tikzpicture}[scale=1.2]
\begin{pgfonlayer}{bg}
	\node[scale=1.2] (00) at (-1,0) {$A$}; 
	\node[scale=1.2] (01) at (0,-1) {$B \lor B$};
	\node[scale=1.2] (10) at (0,1) {$\Sigma^3 (F_s \lor F_s)$};
	\node[scale=1.2] (11) at (1,0) {$IB$};
	\draw[1cinc] (00) to node[auto,swap] {} (01);
	\draw[1cinc] (01) to node[auto,swap] {} (11);
	\draw[1cinc] (00) to node[auto] {} (10);
	\draw[1cinc] (10) to node[auto] {} (11);
	\draw[edge] (0.2,0) to (0.5, -0.3);
	\draw[edge] (0.2,0) to (0.5, 0.3);
\end{pgfonlayer} 
\end{tikzpicture}
\end{equation*}
	in the category of pointed computads.
\end{enumerate}
This does not include the axiom that the dualities induced by the two Frobenius algebras are equal, which is equivalent to the antipode of the bialgebras being the identity. From preliminary results, it appears that to complete the compositional presentation, certain ``skewed'' sliding rules are needed, which can be obtained by working in a cubical, rather than globular, setting, in which different directions are not \emph{a priori} equivalent. We leave this to further work.

\section{Conclusions and outlook}
In this paper, we introduced a basic language for composing higher-dimensional algebraic theories, embodied by computads, in the way that topological spaces can be composed, and demonstrated how simple constructions correspond to common algebraic interactions. 

This language is able to account for the topological differences between the interactions that produce, for instance, Frobenius algebras and bialgebras, respectively, from monoids and comonoids, in a way that earlier compositional frameworks, built in a strictly 2-categorical setting, could not. It seems also remarkably ``inductive'': we showed how to obtain complicated, higher-dimensional coherence diagrams by performing obvious compositions, and then just calculating. This goes in favour of the framework having a heuristic value.

Furthermore, every theory we have constructed in this paper has been obtained from copies of the directed interval $I$ through five basic operations: disjoint union, tensor product, identification of congruent cells, quotient by a subspace, and reversal of cells. By keeping the number of operations contained, we can hope to prove general theorems of the form
\begin{itemize}
	\item[] \emph{If the theory $T$ is obtained from the theories $T_1, \ldots, T_n$, which have the property $P$, by the operation $x$, then $T$ has the property $P'$},
\end{itemize}
and use them to prove interesting facts about interesting theories. By contrast, where cells of arbitrary shape have to be added by hand, as in \cite{bonchi2014hopf, duncan2016interacting}, it is unclear how one could obtain general results.

There are many directions in which to proceed from here. One path is purely incremental: finding more examples, analysing different theories containing different interactions, and trying to describe them in the language of directed topology. As a follow-up to the description of interacting bialgebras, we are particularly interested in a fully topological account of the ZW calculus \cite{hadzihasanovic2015diagrammatic}, including the theory of Hopf algebras.

On a higher level of abstraction, the ideas of this paper could serve as a link between the use of homotopical and homological methods in rewriting theory, and the tools of directed and nonabelian topology developed in monographs such as \cite{grandis2009directed, brown2011nonabelian}. Compositionality is a powerful calculational tool in algebraic topology, from the Seifert-van Kampen Theorem, to the Mayer-Vietoris Theorem, through monoidality of homology functors, so there is a clear potential gain in sight. However, we note that both the latter sources show a clear preference for cubical methods, so a cubical approach to computads might be required to make calculations simpler.

There is, then, the issue of strictness, which does not allow this sort of computads to directly present braidings and other intermediate degrees of commutativity, of the kind that has been extensively studied in the theory of topological operads \cite{markl2007operads}. This could be tackled by resorting to different, weaker notions of computad and higher category; otherwise, in the spirit of \cite{hermida2000representable, hermida2001coherent}, rather than relying on notions of weakness ``from the outside'', it might be conceptually more rewarding to develop them within the simpler combinatorics of strict computads, by imposing various representability conditions.

Finally, as picturing things by diagrams becomes harder and harder in high dimensions, a computational aid may be useful: it could be worth developing an extension of Globular \cite{bar2016globular} to automatise certain compositions.

\section*{Acknowledgments}

The author is supported by an EPSRC Doctoral Training Grant. Thanks to Bob Coecke, Jamie Vicary, Paul-Andr\'e Melli\`es, Samson Abramsky, Stefano Gogioso, Alex Kavvos, and Dominic Verdon for useful discussions and suggestions at various points in the development of these ideas, to Dan Marsden for help with the typesetting of string diagrams, and to the referees for their feedback on a previous draft of this paper, which contained some technical mistakes.

\bibliographystyle{eptcs}
\bibliography{biblinew}

\begin{thebibliography}{10}
\providecommand{\bibitemdeclare}[2]{}
\providecommand{\surnamestart}{}
\providecommand{\surnameend}{}
\providecommand{\urlprefix}{Available at }
\providecommand{\url}[1]{\texttt{#1}}
\providecommand{\href}[2]{\texttt{#2}}
\providecommand{\urlalt}[2]{\href{#1}{#2}}
\providecommand{\doi}[1]{doi:\urlalt{http://dx.doi.org/#1}{#1}}
\providecommand{\bibinfo}[2]{#2}

\bibitemdeclare{article}{al2002multiple}
\bibitem{al2002multiple}
\bibinfo{author}{F.~A. \surnamestart Al-Agl\surnameend},
  \bibinfo{author}{R.~\surnamestart Brown\surnameend} \&
  \bibinfo{author}{R.~\surnamestart Steiner\surnameend} (\bibinfo{year}{2002}):
  \emph{\bibinfo{title}{Multiple categories: the equivalence of a globular and
  a cubical approach}}.
\newblock {\sl \bibinfo{journal}{Advances in Mathematics}}
  \bibinfo{volume}{170}(\bibinfo{number}{1}), pp. \bibinfo{pages}{71--118},
  \doi{10.1006/aima.2001.2069}.

\bibitemdeclare{article}{backens2014zx}
\bibitem{backens2014zx}
\bibinfo{author}{M.~\surnamestart Backens\surnameend} (\bibinfo{year}{2014}):
  \emph{\bibinfo{title}{The ZX-calculus is complete for stabilizer quantum
  mechanics}}.
\newblock {\sl \bibinfo{journal}{New Journal of Physics}}
  \bibinfo{volume}{16}(\bibinfo{number}{9}), p. \bibinfo{pages}{093021},
  \doi{10.1088/1367-2630/16/9/093021}.

\bibitemdeclare{misc}{bar2016globular}
\bibitem{bar2016globular}
\bibinfo{author}{K.~\surnamestart Bar\surnameend},
  \bibinfo{author}{A.~\surnamestart Kissinger\surnameend} \&
  \bibinfo{author}{J.~\surnamestart Vicary\surnameend}:
  \emph{\bibinfo{title}{The \emph{Globular} proof assistant}}.
\newblock \bibinfo{howpublished}{\url{http://ncatlab.org/nlab/show/Globular}}.

\bibitemdeclare{misc}{batanin2002computads}
\bibitem{batanin2002computads}
\bibinfo{author}{M.~\surnamestart Batanin\surnameend} (\bibinfo{year}{2002}):
  \emph{\bibinfo{title}{Computads and slices of operads}}.
\newblock \urlprefix\url{https://arxiv.org/abs/math/0209035}.

\bibitemdeclare{inproceedings}{beck1969distributive}
\bibitem{beck1969distributive}
\bibinfo{author}{J.~\surnamestart Beck\surnameend} (\bibinfo{year}{1969}):
  \emph{\bibinfo{title}{Distributive laws}}.
\newblock In: {\sl \bibinfo{booktitle}{Seminar on triples and categorical
  homology theory}}, \bibinfo{organization}{Springer}, pp.
  \bibinfo{pages}{119--140}, \doi{10.1007/BFb0083084}.

\bibitemdeclare{book}{boardman1973homotopy}
\bibitem{boardman1973homotopy}
\bibinfo{author}{J.~M. \surnamestart Boardman\surnameend} \&
  \bibinfo{author}{R.~M. \surnamestart Vogt\surnameend} (\bibinfo{year}{1973}):
  \emph{\bibinfo{title}{Homotopy invariant algebraic structures on topological
  spaces}}.
\newblock \bibinfo{publisher}{Springer}, \doi{10.1007/978-3-642-54830-7\_23}.

\bibitemdeclare{inproceedings}{bonchi2014interacting}
\bibitem{bonchi2014interacting}
\bibinfo{author}{F.~\surnamestart Bonchi\surnameend},
  \bibinfo{author}{P.~\surnamestart Soboci{\'n}ski\surnameend} \&
  \bibinfo{author}{F.~\surnamestart Zanasi\surnameend} (\bibinfo{year}{2014}):
  \emph{\bibinfo{title}{Interacting Bialgebras Are {F}robenius}}.
\newblock In: {\sl \bibinfo{booktitle}{Foundations of Software Science and
  Computation Structures}}, \bibinfo{organization}{Springer}, pp.
  \bibinfo{pages}{351--365}, \doi{10.1007/978-3-642-54830-7\_23}.

\bibitemdeclare{article}{bonchi2014hopf}
\bibitem{bonchi2014hopf}
\bibinfo{author}{F.~\surnamestart Bonchi\surnameend},
  \bibinfo{author}{P.~\surnamestart Sobocinski\surnameend} \&
  \bibinfo{author}{F.~\surnamestart Zanasi\surnameend} (\bibinfo{year}{2017}):
  \emph{\bibinfo{title}{Interacting Hopf algebras}}.
\newblock {\sl \bibinfo{journal}{Journal of Pure and Applied Algebra}}
  \bibinfo{volume}{221}(\bibinfo{number}{1}), pp. \bibinfo{pages}{144 -- 184},
  \doi{10.1016/j.jpaa.2016.06.002}.

\bibitemdeclare{article}{brown1981algebra}
\bibitem{brown1981algebra}
\bibinfo{author}{R.~\surnamestart Brown\surnameend} \& \bibinfo{author}{P.~J.
  \surnamestart Higgins\surnameend} (\bibinfo{year}{1981}):
  \emph{\bibinfo{title}{On the algebra of cubes}}.
\newblock {\sl \bibinfo{journal}{Journal of Pure and Applied Algebra}}
  \bibinfo{volume}{21}(\bibinfo{number}{3}), pp. \bibinfo{pages}{233--260},
  \doi{10.1016/0022-4049(81)90018-9}.

\bibitemdeclare{article}{brown1987tensor}
\bibitem{brown1987tensor}
\bibinfo{author}{R.~\surnamestart Brown\surnameend} \& \bibinfo{author}{P.~J.
  \surnamestart Higgins\surnameend} (\bibinfo{year}{1987}):
  \emph{\bibinfo{title}{Tensor products and homotopies for $\omega$-groupoids
  and crossed complexes}}.
\newblock {\sl \bibinfo{journal}{Journal of Pure and Applied Algebra}}
  \bibinfo{volume}{47}(\bibinfo{number}{1}), pp. \bibinfo{pages}{1--33},
  \doi{10.1016/0022-4049(87)90099-5}.

\bibitemdeclare{book}{brown2011nonabelian}
\bibitem{brown2011nonabelian}
\bibinfo{author}{R.~\surnamestart Brown\surnameend}, \bibinfo{author}{P.~J.
  \surnamestart Higgins\surnameend} \& \bibinfo{author}{R.~\surnamestart
  Sivera\surnameend} (\bibinfo{year}{2011}): \emph{\bibinfo{title}{Nonabelian
  algebraic topology}}.
\newblock \bibinfo{publisher}{European Mathematical Society},
  \doi{10.4171/083}.

\bibitemdeclare{misc}{buckley2015orientals}
\bibitem{buckley2015orientals}
\bibinfo{author}{M.~\surnamestart Buckley\surnameend} \&
  \bibinfo{author}{R.~\surnamestart Garner\surnameend} (\bibinfo{year}{2015}):
  \emph{\bibinfo{title}{Orientals and cubes, inductively}}.
\newblock \urlprefix\url{https://arxiv.org/abs/1509.00618}.

\bibitemdeclare{article}{burroni1993higher}
\bibitem{burroni1993higher}
\bibinfo{author}{A.~\surnamestart Burroni\surnameend} (\bibinfo{year}{1993}):
  \emph{\bibinfo{title}{Higher-dimensional word problems with applications to
  equational logic}}.
\newblock {\sl \bibinfo{journal}{Theoretical Computer Science}}
  \bibinfo{volume}{115}(\bibinfo{number}{1}), pp. \bibinfo{pages}{43--62},
  \doi{10.1016/0304-3975(93)90054-W}.

\bibitemdeclare{incollection}{coecke2008interacting}
\bibitem{coecke2008interacting}
\bibinfo{author}{B.~\surnamestart Coecke\surnameend} \&
  \bibinfo{author}{R.~\surnamestart Duncan\surnameend} (\bibinfo{year}{2008}):
  \emph{\bibinfo{title}{Interacting quantum observables}}.
\newblock In: {\sl \bibinfo{booktitle}{Automata, Languages and Programming}},
  \bibinfo{publisher}{Springer}, pp. \bibinfo{pages}{298--310},
  \doi{10.1007/978-3-540-70583-3\_25}.

\bibitemdeclare{book}{coecke2015picturing}
\bibitem{coecke2015picturing}
\bibinfo{author}{B.~\surnamestart Coecke\surnameend} \&
  \bibinfo{author}{A.~\surnamestart Kissinger\surnameend}
  (\bibinfo{year}{2017}): \emph{\bibinfo{title}{Picturing quantum processes}}.
\newblock \bibinfo{publisher}{Cambridge University Press}.
\newblock \bibinfo{note}{To appear}.

\bibitemdeclare{article}{crans1995pasting}
\bibitem{crans1995pasting}
\bibinfo{author}{S.~E. \surnamestart Crans\surnameend} (\bibinfo{year}{1995}):
  \emph{\bibinfo{title}{Pasting schemes for the monoidal biclosed structure on
  $\omega$-Cat}}.
\newblock \bibinfo{note}{Utrecht University}.

\bibitemdeclare{inproceedings}{duncan2016interacting}
\bibitem{duncan2016interacting}
\bibinfo{author}{R.~\surnamestart Duncan\surnameend} \&
  \bibinfo{author}{K.~\surnamestart Dunne\surnameend} (\bibinfo{year}{2016}):
  \emph{\bibinfo{title}{Interacting Frobenius Algebras Are Hopf}}.
\newblock In: {\sl \bibinfo{booktitle}{Proceedings of the 31st Annual ACM/IEEE
  Symposium on Logic in Computer Science}}, \bibinfo{series}{LICS '16}, pp.
  \bibinfo{pages}{535--544}, \doi{10.1145/2933575.2934550}.

\bibitemdeclare{misc}{dunn2016surface}
\bibitem{dunn2016surface}
\bibinfo{author}{L.~\surnamestart Dunn\surnameend} \&
  \bibinfo{author}{J.~\surnamestart Vicary\surnameend} (\bibinfo{year}{2016}):
  \emph{\bibinfo{title}{Coherence for Frobenius pseudomonoids and the geometry
  of linear proofs}}.
\newblock \urlprefix\url{https://arxiv.org/abs/1601.05372v3}.

\bibitemdeclare{article}{eckmann1962group}
\bibitem{eckmann1962group}
\bibinfo{author}{B.~\surnamestart Eckmann\surnameend} \& \bibinfo{author}{P.~J.
  \surnamestart Hilton\surnameend} (\bibinfo{year}{1962}):
  \emph{\bibinfo{title}{Group-like structures in general categories {I}:
  multiplications and comultiplications}}.
\newblock {\sl \bibinfo{journal}{Mathematische Annalen}}
  \bibinfo{volume}{145}(\bibinfo{number}{3}), pp. \bibinfo{pages}{227--255},
  \doi{10.1007/BF01451367}.

\bibitemdeclare{book}{grandis2009directed}
\bibitem{grandis2009directed}
\bibinfo{author}{M.~\surnamestart Grandis\surnameend} (\bibinfo{year}{2009}):
  \emph{\bibinfo{title}{Directed Algebraic Topology: Models of non-reversible
  worlds}}.
\newblock \bibinfo{volume}{13}, \bibinfo{publisher}{Cambridge University
  Press}, \doi{10.1017/CBO9780511657474}.

\bibitemdeclare{article}{guiraud2006three}
\bibitem{guiraud2006three}
\bibinfo{author}{Y.~\surnamestart Guiraud\surnameend} (\bibinfo{year}{2006}):
  \emph{\bibinfo{title}{The three dimensions of proofs}}.
\newblock {\sl \bibinfo{journal}{Annals of Pure and Applied Logic}}
  \bibinfo{volume}{141}(\bibinfo{number}{1}), pp. \bibinfo{pages}{266--295},
  \doi{10.1016/j.apal.2005.12.012}.

\bibitemdeclare{inproceedings}{hadzihasanovic2015diagrammatic}
\bibitem{hadzihasanovic2015diagrammatic}
\bibinfo{author}{A.~\surnamestart Hadzihasanovic\surnameend}
  (\bibinfo{year}{2015}): \emph{\bibinfo{title}{A Diagrammatic Axiomatisation
  for Qubit Entanglement}}.
\newblock In: {\sl \bibinfo{booktitle}{Logic in Computer Science (LICS), 2015
  30th Annual ACM/IEEE Symposium on}}, \bibinfo{organization}{IEEE}, pp.
  \bibinfo{pages}{573--584}, \doi{10.1109/LICS.2015.59}.

\bibitemdeclare{article}{hermida2000representable}
\bibitem{hermida2000representable}
\bibinfo{author}{C.~\surnamestart Hermida\surnameend} (\bibinfo{year}{2000}):
  \emph{\bibinfo{title}{Representable multicategories}}.
\newblock {\sl \bibinfo{journal}{Advances in Mathematics}}
  \bibinfo{volume}{151}(\bibinfo{number}{2}), pp. \bibinfo{pages}{164--225},
  \doi{10.1006/aima.1999.1877}.

\bibitemdeclare{article}{hermida2001coherent}
\bibitem{hermida2001coherent}
\bibinfo{author}{C.~\surnamestart Hermida\surnameend} (\bibinfo{year}{2001}):
  \emph{\bibinfo{title}{From coherent structures to universal properties}}.
\newblock {\sl \bibinfo{journal}{Journal of Pure and Applied Algebra}}
  \bibinfo{volume}{165}(\bibinfo{number}{1}), pp. \bibinfo{pages}{7--61},
  \doi{10.1016/S0022-4049(01)00008-1}.

\bibitemdeclare{article}{heunen2012lectures}
\bibitem{heunen2012lectures}
\bibinfo{author}{C.~\surnamestart Heunen\surnameend} \&
  \bibinfo{author}{J.~\surnamestart Vicary\surnameend} (\bibinfo{year}{2012}):
  \emph{\bibinfo{title}{Lectures on categorical quantum mechanics}}.
\newblock {\sl \bibinfo{journal}{Computer Science Department. Oxford
  University}}.

\bibitemdeclare{article}{hinze2016equational}
\bibitem{hinze2016equational}
\bibinfo{author}{R.~\surnamestart Hinze\surnameend} \&
  \bibinfo{author}{D.~\surnamestart Marsden\surnameend} (\bibinfo{year}{2016}):
  \emph{\bibinfo{title}{Equational reasoning with lollipops, forks, cups, caps,
  snakes, and speedometers}}.
\newblock {\sl \bibinfo{journal}{Journal of Logical and Algebraic Methods in
  Programming}}, \doi{10.1016/j.jlamp.2015.12.004}.

\bibitemdeclare{article}{lack2004composing}
\bibitem{lack2004composing}
\bibinfo{author}{S.~\surnamestart Lack\surnameend} (\bibinfo{year}{2004}):
  \emph{\bibinfo{title}{Composing {PROP}s}}.
\newblock {\sl \bibinfo{journal}{Theory and Applications of Categories}}
  \bibinfo{volume}{13}(\bibinfo{number}{9}), pp. \bibinfo{pages}{147--163}.

\bibitemdeclare{article}{lafont2007algebra}
\bibitem{lafont2007algebra}
\bibinfo{author}{Y.~\surnamestart Lafont\surnameend} (\bibinfo{year}{2007}):
  \emph{\bibinfo{title}{Algebra and geometry of rewriting}}.
\newblock {\sl \bibinfo{journal}{Applied Categorical Structures}}
  \bibinfo{volume}{15}(\bibinfo{number}{4}), pp. \bibinfo{pages}{415--437},
  \doi{10.1007/s10485-007-9083-6}.

\bibitemdeclare{article}{lafont2009polygraphic}
\bibitem{lafont2009polygraphic}
\bibinfo{author}{Y.~\surnamestart Lafont\surnameend} \&
  \bibinfo{author}{F.~\surnamestart M{\'e}tayer\surnameend}
  (\bibinfo{year}{2009}): \emph{\bibinfo{title}{Polygraphic resolutions and
  homology of monoids}}.
\newblock {\sl \bibinfo{journal}{Journal of Pure and Applied Algebra}}
  \bibinfo{volume}{213}(\bibinfo{number}{6}), pp. \bibinfo{pages}{947--968},
  \doi{10.1016/j.jpaa.2008.10.005}.

\bibitemdeclare{article}{maclane1963natural}
\bibitem{maclane1963natural}
\bibinfo{author}{S.~\surnamestart MacLane\surnameend} (\bibinfo{year}{1963}):
  \emph{\bibinfo{title}{Natural associativity and commutativity}}.
\newblock {\sl \bibinfo{journal}{Rice Institute Pamphlet-Rice University
  Studies}} \bibinfo{volume}{49}(\bibinfo{number}{4}).

\bibitemdeclare{misc}{makkai2001comparing}
\bibitem{makkai2001comparing}
\bibinfo{author}{M.~\surnamestart Makkai\surnameend} (\bibinfo{year}{2005}):
  \emph{\bibinfo{title}{The word problem for computads}}.
\newblock \bibinfo{note}{Available on the author's web page
  \url{http://www.math.mcgill.ca/makkai/}}.

\bibitemdeclare{article}{markl2008operads}
\bibitem{markl2008operads}
\bibinfo{author}{M.~\surnamestart Markl\surnameend} (\bibinfo{year}{2008}):
  \emph{\bibinfo{title}{Operads and PROPS}}.
\newblock {\sl \bibinfo{journal}{Handbook of Algebra}} \bibinfo{volume}{5}, pp.
  \bibinfo{pages}{87--140}, \doi{10.1016/S1570-7954(07)05002-4}.

\bibitemdeclare{book}{markl2007operads}
\bibitem{markl2007operads}
\bibinfo{author}{M.~\surnamestart Markl\surnameend},
  \bibinfo{author}{S.~\surnamestart Shnider\surnameend} \&
  \bibinfo{author}{J.~D. \surnamestart Stasheff\surnameend}
  (\bibinfo{year}{2007}): \emph{\bibinfo{title}{Operads in algebra, topology
  and physics}}.
\newblock \bibinfo{volume}{96}, \bibinfo{publisher}{American Mathematical
  Soc.}, \doi{10.1090/surv/096}.

\bibitemdeclare{article}{metayer2003resolutions}
\bibitem{metayer2003resolutions}
\bibinfo{author}{F.~\surnamestart M{\'e}tayer\surnameend}
  (\bibinfo{year}{2003}): \emph{\bibinfo{title}{Resolutions by polygraphs}}.
\newblock {\sl \bibinfo{journal}{Theory and Applications of Categories}}
  \bibinfo{volume}{11}(\bibinfo{number}{7}), pp. \bibinfo{pages}{148--184}.

\bibitemdeclare{misc}{mimram2014towards}
\bibitem{mimram2014towards}
\bibinfo{author}{S.~\surnamestart Mimram\surnameend} (\bibinfo{year}{2014}):
  \emph{\bibinfo{title}{Towards 3-dimensional rewriting theory}}.
\newblock \urlprefix\url{https://arxiv.org/abs/1410.2901}.

\bibitemdeclare{misc}{perk2006yang}
\bibitem{perk2006yang}
\bibinfo{author}{J.H.H. \surnamestart Perk\surnameend} \&
  \bibinfo{author}{H.~\surnamestart Au-Yang\surnameend} (\bibinfo{year}{2006}):
  \emph{\bibinfo{title}{Yang-{B}axter equations}}.
\newblock \urlprefix\url{https://arxiv.org/abs/math-ph/0606053}.

\bibitemdeclare{article}{selinger2011finite}
\bibitem{selinger2011finite}
\bibinfo{author}{P.~\surnamestart Selinger\surnameend} (\bibinfo{year}{2011}):
  \emph{\bibinfo{title}{Finite dimensional Hilbert spaces are complete for
  dagger compact closed categories}}.
\newblock {\sl \bibinfo{journal}{Electronic Notes in Theoretical Computer
  Science}} \bibinfo{volume}{270}(\bibinfo{number}{1}), pp.
  \bibinfo{pages}{113--119}, \doi{10.1016/j.entcs.2011.01.010}.

\bibitemdeclare{incollection}{selinger2011survey}
\bibitem{selinger2011survey}
\bibinfo{author}{P.~\surnamestart Selinger\surnameend} (\bibinfo{year}{2011}):
  \emph{\bibinfo{title}{A survey of graphical languages for monoidal
  categories}}.
\newblock In: {\sl \bibinfo{booktitle}{New structures for physics}},
  \bibinfo{publisher}{Springer}, pp. \bibinfo{pages}{289--355}.

\bibitemdeclare{article}{steiner2004omega}
\bibitem{steiner2004omega}
\bibinfo{author}{R.~\surnamestart Steiner\surnameend} (\bibinfo{year}{2004}):
  \emph{\bibinfo{title}{Omega-categories and chain complexes}}.
\newblock {\sl \bibinfo{journal}{Homology, Homotopy and Applications}}
  \bibinfo{volume}{6}(\bibinfo{number}{1}), pp. \bibinfo{pages}{175--200},
  \doi{10.4310/HHA.2004.v6.n1.a12}.

\bibitemdeclare{article}{street1976limits}
\bibitem{street1976limits}
\bibinfo{author}{R.~\surnamestart Street\surnameend} (\bibinfo{year}{1976}):
  \emph{\bibinfo{title}{Limits indexed by category-valued 2-functors}}.
\newblock {\sl \bibinfo{journal}{Journal of Pure and Applied Algebra}}
  \bibinfo{volume}{8}(\bibinfo{number}{2}), pp. \bibinfo{pages}{149--181},
  \doi{10.1016/0022-4049(76)90013-X}.

\bibitemdeclare{article}{street1987algebra}
\bibitem{street1987algebra}
\bibinfo{author}{R.~\surnamestart Street\surnameend} (\bibinfo{year}{1987}):
  \emph{\bibinfo{title}{The algebra of oriented simplexes}}.
\newblock {\sl \bibinfo{journal}{Journal of Pure and Applied Algebra}}
  \bibinfo{volume}{49}(\bibinfo{number}{3}), pp. \bibinfo{pages}{283--335},
  \doi{10.1016/0022-4049(87)90137-X}.

\bibitemdeclare{misc}{tubella2016subatomic}
\bibitem{tubella2016subatomic}
\bibinfo{author}{A.~A. \surnamestart Tubella\surnameend} \&
  \bibinfo{author}{A.~\surnamestart Guglielmi\surnameend}
  (\bibinfo{year}{2016}): \emph{\bibinfo{title}{Subatomic Proof Systems}}.
\newblock \bibinfo{note}{Available on the author's web page
  \url{http://alessio.guglielmi.name/res/cos/}}.

\bibitemdeclare{article}{weiss2011operads}
\bibitem{weiss2011operads}
\bibinfo{author}{I.~\surnamestart Weiss\surnameend} (\bibinfo{year}{2011}):
  \emph{\bibinfo{title}{From operads to dendroidal sets}}.
\newblock {\sl \bibinfo{journal}{Mathematical foundations of quantum field
  theory and perturbative string theory}} \bibinfo{volume}{83}, pp.
  \bibinfo{pages}{31--70}, \doi{10.1090/pspum/083/2742425}.

\bibitemdeclare{misc}{zeng2015abstract}
\bibitem{zeng2015abstract}
\bibinfo{author}{W.~\surnamestart Zeng\surnameend} (\bibinfo{year}{2015}):
  \emph{\bibinfo{title}{The Abstract Structure of Quantum Algorithms}}.
\newblock \urlprefix\url{https://arxiv.org/abs/1512.08062}.

\end{thebibliography}
\begin{appendices}
\section{Strict $\omega$-categories} \label{sec:globular}

We recall the definition of strict, globular $\omega$-category.

\begin{dfn}
A \emph{strict $\omega$-category} $\mathcal{C}$ is a set together with unary border operators $\partial^+_n$ (output $n$-border), $\partial^-_n$ (input $n$-border), and partial binary compositions $\cp{n}$, for all $n \geq 0$, satisfying the following axioms:
\begin{enumerate}
	\item for all $\sigma, \tau$ in $\mathcal{C}$, $\sigma \cp{n} \tau$ is defined if and only if $\partial^+_n\sigma = \partial^-_n\tau$;
	\item whenever both sides are defined, and $n \neq m$,
		\begin{align*}
			& (\sigma \cp{n} \tau) \cp{n} \rho = \sigma \cp{n} (\tau \cp{n} \rho) & \text{(associativity),}\\ 
			& \sigma \cp{n} \partial^+_n\sigma = \sigma = \partial^-_n\sigma \cp{n} \sigma & \text{(unitality),}\\
			& (\sigma_1 \cp{n} \sigma_2) \cp{m} (\tau_1 \cp{n} \tau_2) = (\sigma_1 \cp{m} \tau_1) \cp{n} (\sigma_2 \cp{m} \tau_2) & \text{(interchange);}
		\end{align*}
	\item for all $n, m \geq 0$, and $\alpha, \beta \in \{+,-\}$,
		\begin{equation*}
			\partial^\beta_m\partial^\alpha_n = \begin{cases}
				\partial^\beta_m\;, & m < n\;, \\
				\partial^\alpha_n\;, & m \geq n\;;
			\end{cases}
		\end{equation*} 
	\item whenever $\sigma \cp{n} \tau$ is defined, and $m \neq n$,
		\begin{align*}
			\partial^-_n(\sigma \cp{n} \tau) & = \partial^-_n\sigma\;, \\
			\partial^+_n(\sigma \cp{n} \tau) & = \partial^+_n\tau\;, \\
			\partial^\alpha_m(\sigma \cp{n} \tau) & = \partial^\alpha_m\sigma \cp{n} \partial^\alpha_m\tau\;;
		\end{align*}
	\item for all $\sigma$ in $\mathcal{C}$, there is a smallest $n$, the \emph{dimension} $d(\sigma)$ of $\sigma$, such that for all $m \geq n$
		\begin{equation*}
			\partial^-_m \sigma = \sigma = \partial^+_m \sigma\;.
		\end{equation*}
\end{enumerate}
Elements of $\mathcal{C}$ are called \emph{cells}; a cell of dimension $n$ is an $n$-cell. For any $\omega$-category $\mathcal{C}$, and $n \geq 0$, the \emph{$n$-skeleton} $\mathcal{C}_n$ of $\mathcal{C}$ is the restriction of $\mathcal{C}$ to cells of dimension $d \leq n$. 

Given two $\omega$-categories $\mathcal{C}$, $\mathcal{D}$, a \emph{functor} $f: \mathcal{C} \to \mathcal{D}$ is a function commuting with border operators and compositions. Functors and $\omega$-categories form a category $\omega\mathbf{Cat}$.
\end{dfn}

\section{Borders in low dimensions} \label{sec:borders}

Let $X$, $Y$ be two computads, $\sigma \in |X|$, $\tau \in |Y|$. We give explicit expressions for the low-dimensional borders of $\sigma \otimes \tau$, in terms of the borders of $\sigma$ and $\tau$.
\begin{align*}
	\bord{0}{-}(\sigma \otimes \tau) \; = & \;\; \bord{0}{-}\sigma \otimes \bord{0}{-}\tau \\
	\bord{0}{+}(\sigma \otimes \tau) \; = & \;\; \bord{0}{+}\sigma \otimes \bord{0}{+}\tau \\
	\bord{1}{-}(\sigma \otimes \tau) \; = & \;\; (\bord{0}{-}\sigma \otimes \bord{1}{-}\tau) \cp{1} (\bord{1}{-}\sigma \otimes \bord{0}{+}\tau) \\
	\bord{1}{+}(\sigma \otimes \tau) \; = & \;\; (\bord{1}{+}\sigma \otimes \bord{0}{-}\tau) \cp{1} (\bord{0}{+}\sigma \otimes \bord{1}{+}\tau) \; 	\\
	\bord{2}{-}(\sigma \otimes \tau) \; = & \;\; \big((\bord{0}{-}\sigma \otimes \bord{2}{-}\tau) \cp{1} (\bord{1}{-}\sigma \otimes \bord{0}{+}\tau) \big)
		\cp{2} ( \bord{1}{-}\sigma \otimes \bord{1}{+}\tau ) 
		\cp{2} \big(( \bord{2}{-}\sigma \otimes \bord{0}{-}\tau) \cp{1} (\bord{0}{+}\sigma \otimes \bord{1}{+}\tau) \big) \\
	\bord{2}{+}(\sigma \otimes \tau) \; = & \;\; \big((\bord{0}{-}\sigma \otimes \bord{1}{-}\tau) \cp{1} (\bord{2}{+}\sigma \otimes \bord{0}{+}\tau) \big)
		\cp{2} ( \bord{1}{+}\sigma \otimes \bord{1}{-}\tau ) 
		\cp{2} \big(( \bord{1}{+}\sigma \otimes \bord{0}{-}\tau) \cp{1} (\bord{0}{+}\sigma \otimes \bord{2}{+}\tau) \big) \\
	\bord{3}{-}(\sigma \otimes \tau) \; = & \;\; \Big(\big((\bord{0}{-}\sigma \otimes \bord{3}{-}\tau) \cp{1} (\bord{1}{-}\sigma \otimes \bord{0}{+}\tau) \big)
		\cp{2} ( \bord{1}{-}\sigma \otimes \bord{1}{+}\tau ) 
		\cp{2} \big(( \bord{2}{-}\sigma \otimes \bord{0}{-}\tau) \cp{1} (\bord{0}{+}\sigma \otimes \bord{1}{+}\tau) \big)\Big) \\
		& \cp{3} \Big((\bord{1}{-}\sigma \otimes \bord{2}{+}\tau) \cp{2} \big(( \bord{2}{-}\sigma \otimes \bord{0}{-}\tau) \cp{1} (\bord{0}{+}\sigma \otimes \bord{1}{+}\tau) \big) \Big) \\
		& \cp{3} \Big((\bord{2}{-}\sigma \otimes \bord{1}{-}\tau) \cp{2} \big(( \bord{1}{+}\sigma \otimes \bord{0}{-}\tau) \cp{1} (\bord{0}{+}\sigma \otimes \bord{2}{+}\tau) \big) \Big) \\
		& \cp{3} \Big(\big((\bord{0}{-}\sigma \otimes \bord{1}{-}\tau) \cp{1} (\bord{3}{-}\sigma \otimes \bord{0}{+}\tau) \big)
		\cp{2} ( \bord{1}{+}\sigma \otimes \bord{1}{-}\tau ) 
		\cp{2} \big(( \bord{1}{+}\sigma \otimes \bord{0}{-}\tau) \cp{1} (\bord{0}{+}\sigma \otimes \bord{2}{+}\tau) \big)\Big) \\
	\bord{3}{+}(\sigma \otimes \tau) \; = & \;\; \Big(\big((\bord{0}{-}\sigma \otimes \bord{2}{-}\tau) \cp{1} (\bord{1}{-}\sigma \otimes \bord{0}{+}\tau) \big)
		\cp{2} ( \bord{1}{-}\sigma \otimes \bord{1}{+}\tau ) 
		\cp{2} \big(( \bord{3}{+}\sigma \otimes \bord{0}{-}\tau) \cp{1} (\bord{0}{+}\sigma \otimes \bord{1}{+}\tau) \big)\Big) \\
		& \cp{3} \Big( \big(( \bord{0}{-}\sigma \otimes \bord{2}{-}\tau) \cp{1} (\bord{1}{-}\sigma \otimes \bord{0}{+}\tau) \big) \cp{2} (\bord{2}{+}\sigma \otimes \bord{1}{+}\tau) \Big) \\
		& \cp{3} \Big( \big(( \bord{0}{-}\sigma \otimes \bord{1}{-}\tau) \cp{1} (\bord{2}{+}\sigma \otimes \bord{0}{+}\tau) \big) \cp{2} (\bord{1}{+}\sigma \otimes \bord{2}{-}\tau) \Big) \\
		& \cp{3} \Big(\big((\bord{0}{-}\sigma \otimes \bord{1}{-}\tau) \cp{1} (\bord{2}{+}\sigma \otimes \bord{0}{+}\tau) \big)
		\cp{2} ( \bord{1}{+}\sigma \otimes \bord{1}{-}\tau ) 
		\cp{2} \big(( \bord{1}{+}\sigma \otimes \bord{0}{-}\tau) \cp{1} (\bord{0}{+}\sigma \otimes \bord{3}{+}\tau) \big)\Big) \\		
\end{align*}

\end{appendices}

\end{document}